\documentclass[a4paper]{amsart}

\usepackage{lipsum}
\usepackage{amsfonts}
\usepackage{graphicx}
\usepackage{epstopdf}
\usepackage{algorithmic}
\ifpdf
  \DeclareGraphicsExtensions{.eps,.pdf,.png,.jpg}
\else
  \DeclareGraphicsExtensions{.eps}
\fi

\title[MACROSCOPIC MULTILANE MODELS FOR TRAFFIC FLOW]{Derivation and stability analysis of  macroscopic multi-lane models for vehicular traffic flow}

\author[M.Piu]{Matteo Piu}
\address{Corresponding author. Department of Basic and Applied Sciences for Engineering, Sapienza - University of Rome, Rome, Italy}
\email{matteo.piu@uniroma1.it}
\author[M.Herty]{Michael Herty}
\address{Institute of Applied Mathematics, RWTH Aachen University, Aachen, Germany}
\email{herty@igpm.rwth-aachen.de}
\author[G.Puppo]{Gabriella Puppo}
\address{Department of Mathematics, Sapienza - University of Rome, Rome, Italy}
\email{gabriella.puppo@uniroma1.it}

\usepackage{amsopn}

\usepackage{multirow}
\newcommand{\patP}{crosshatch}
\newcommand{\patUnoMenoPmezzi}{crosshatch dots}
\usepackage{tikz,pgfplots}
\usetikzlibrary{intersections,calc,patterns}
\newcommand{\de}{\partial}

\newcommand{\sign}{\text{sign}}
\newcommand{\epi}{\varepsilon}
\usepackage{amssymb}
\allowdisplaybreaks

\usepackage{xcolor,subfigure}
\pgfplotsset{compat=1.17} 
\newtheorem{theorem}{Theorem}[section]

\theoremstyle{definition}
\newtheorem{definition}[theorem]{Definition}
\newtheorem{remark}{Remark}

\usepackage{float}

\usepackage[style=numeric,backend=biber,maxnames=10,doi=false,isbn=false,url=false,sortcites=true]{biblatex}
\begin{document}

\maketitle

\begin{abstract}
The mathematical modeling and the stability analysis of multi-lane traffic in the macroscopic scale is considered. We propose a new first order model derived from microscopic dynamics with lane changing, leading to a coupled system of hyperbolic balance laws. The macroscopic limit is derived without assuming ad hoc space and time scalings. The analysis of the stability of the equilibria of the model is discussed. The proposed numerical tests confirm the theoretical findings between the macroscopic and  microscopic modeling, and the results of the stability analysis.
\end{abstract}





\section{Introduction}

Interest in the modeling of dynamics of traffic flow dates back to the first half of the twentieth century and the related mathematical literature is today extensive, see e.g.~\cite{MR1600250,Helbing20011067}. Three different scales of description of  vehicular traffic are typically considered, leading to microscopic, mesoscopic and macroscopic mathematical models; see e.g.~\cite{piccoli2009ENCYCLOPEDIA,albi2019vehicular} 

In this paper we focus on the macroscopic scale, inspired by fluid dynamics equations. Here, traffic flow is represented by the study of average quantities, such as the traffic density. First macroscopic models for traffic flow for single lane roads were proposed in the 1950s by Lighthill and Whitham~\cite{lighthill1955PRSL} and  Richards~\cite{richards1956OR}. This model is derived by the mass conservation principle and it is described by a scalar conservation law of the form
\begin{equation} \label{eq:lwr}
\de_t \rho + \de_x f(\rho) = 0,
\end{equation}
where $\rho=\rho(x,t)$ represents the density of vehicles, with $\rho \in [0,\rho^{\max}]$, and $\rho^{\max}$ is the maximum density of vehicles on the road; $(x,t)\in \mathbb{R}\times \mathbb{R}^+$ represents the space and time variables, respectively. The nonlinear flux function $f(\rho)=\rho v(\rho)$ is a given function of the density and it is called the fundamental diagram, where $v$ is the average velocity of the flow, and, typically, it is a given function of the density. Equation~\eqref{lwr} is endowed with an initial condition $\rho(x,0)=\rho_0(x)$. Improvements and further evolution of the basic macroscopic single lane description~\eqref{eq:lwr} have been proposed over the years by several authors, e.g.~Aw and Rascle~\cite{aw2000SIAP} and, independently, by Zhang~\cite{Zhang2002}. In second order models the scalar conservation law~\eqref{eq:lwr} is coupled to an additional equation for the evolution of the average speed $v$.

Single lane models are scarcely predictive for multi-lane roads, such as highways and freeways. Therefore, mathematical models have been extended in order to describe vehicular dynamics on multiple lanes. This allows us to analyze  effects of lane changing maneuvers on the global behavior and the evolution of traffic flow, such as the appearance of unstable phenomena. Typically, multi-lane models are based on systems of one-dimensional equations for the density of each lane and describe lane changes through source terms modeled directly at the macroscopic level, see, e.g., \cite{GoatinRossi2019,HoldenRisebro2019,SongKarni2019}, using empirical considerations.
In these models the source term is assumed to be proportional to local density on both the current and the target lane. Other fluid dynamics models, instead, describe the cumulative density on all lanes using a two-dimensional system of conservation
laws~\cite{SukhinovaTrapeznikovaChetverushkinChurbanova2009,HertyMoutariVisconti2018}. In
those studies, the modeling of traffic dynamics assumes that vehicles move to lanes with faster speed or lower density, and that the evolution for the lateral velocity is proportional to the local density and the mean speed along the road.

In this paper we consider a different approach. We present a novel first order macroscopic model for multi-lane traffic flow whose derivation is motivated by a microscopic scale description of lane changing dynamics. The position of vehicles is modeled by a system of first order ordinary differential equations, namely by continuous dynamics, whereas lane changes are formulated as discrete events and occur instantaneously, contrary to other approaches based on the definition of a cool-down time~\cite{GongPiccoliVisconti,ChiriGongPiccoli}. The presence of both continuous and discrete dynamics gives us a hybrid system, e.g.~see~\cite{garavello2005hybrid,4806347,piccoli1998hybrid}. The macroscopic model is described by a system of hyperbolic balance laws for the evolution of the density in each lane, including source terms. Here, lane changes are obtained as a macroscopic limit of microscopic lane changing maneuvers. The microscopic-to-macroscopic limit  does not need to postulate specific ad hoc scalings \cite{aw2002SIAP}. In particular we study the effects of  lane changing from the point of view of the affected lane, deriving specific source terms for the evolutionary density equation. In this sense, let us suppose that a vehicle performs  a lane change, and we imagine  observing this occurrence from the perspective of the vehicle in the new lane, immediately behind the vehicle that has just changed. Inevitably, the local density seen by such a vehicle will increase and the derived source term, with respect to the macroscopic variables, aims to model the increase of the density in the lane affected by the lane change in terms of the  density before the change in the same lane. In this way we obtain a source term describing the evolution of the density in a macroscopic scale. The model is analyzed by classifying all possible equilibrium solutions for the two-lane case and the stability with respect to small perturbations.

The paper is organized as follows. Section~\ref{sec:micro} is devoted to the definition of the microscopic framework, with a description of the lane changing conditions, based on safety and incentive criteria as, e.g., in~\cite{PiuPuppo2022,GongPiccoliVisconti}. Then, in Section~\ref{sec:derivation} we derive the first order continuum model as a macroscopic limit of the microscopic continuous dynamics and of the discrete events describing the lane changing. Section~\ref{sec:analysis}  presents the characterization of the equilibria of the macroscopic model in the case of a two-lane road and the analysis of the stability under some specific perturbations. Numerical simulations in Section~\ref{sec:numerics} show  evidence of the theoretical findings, as well as the consistency between the microscopic and the macroscopic models, and an application to the lane closure problem. Finally, in Section~\ref{sec:conclusion} we summarize the results. 

\section{First order multi-lane microscopic model} \label{sec:micro}

In this section we present a first order microscopic multi-lane model characterized by continuous dynamics combined with discrete events for  lane changing.

\subsection{Microscopic dynamics with lane changing}

In the following we consider a homogeneous population of $N\in\mathbb{N}$ vehicles and we denote by $x_n=x_n(t)\in\mathbb{R}$ the time-dependent position of the $n$-th vehicle with $n=1,\dots,N$. More precisely, we denote by $x_n$  the position of the rear bumper, therefore the position of the front bumper is given by $x_n+l$, where $l>0$ denotes the length of the vehicle which is assumed constant for all vehicles.

\begin{figure}[t]
\centering
\begin{tikzpicture}[scale=0.8]
[auto,
 block/.style ={rectangle, draw=blue, thick, fill=blue!20},
 block1/.style ={rectangle, draw=green, thick, fill=green!20},
 block2/.style ={rectangle, draw=red, thick, fill=red!20}
]
\draw [dashed](0,3)--(8,3);
\draw [dashed](0,2)--(8,2);
\draw [dashed](0,1)--(8,1);
\draw [dashed](0,0)--(8,0);
\draw (0,-1)--(8,-1);
\draw (0,4)--(8,4);

\path(-1,0.5)node{lane $j-1$};
\path (-1,1.5)node{lane $j$};
\path (-1,2.5)node{lane $j+1$};

\draw (1,2.1) rectangle (2.5,2.9);
\draw (1+1.5/2,2.5)node{$p_n^{j+1}$};

\draw (6,2.1) rectangle (7.5,2.9);
\draw (6+1.5/2,2.5)node{$s_n^{j+1}$};

\draw (0.5,1.1) rectangle (2,1.9);
\draw (0.5+1.5/2,1.5)node{$p_n^{j}$};

\draw (3,1.1) rectangle (4.5,1.9);
\draw (3+1.5/2,1.5)node{$n$};

\draw (5.6,1.1) rectangle (5.6+1.5,1.9);
\draw (5.6+1.5/2,1.5)node{$s_n^{j}$};

\draw (1.3,0.1) rectangle (2.8,0.9);
\draw (1.3+1.5/2,0.5)node{$p_n^{j-1}$};

\draw (5.2,0.1) rectangle (5.2+1.5,0.9);
\draw (5.2+1.5/2,0.5)node{$s_n^{j-1}$};

\end{tikzpicture}
\caption{Schematic representation of a multi-lane road. Here, the reference vehicle is $n$, travelling in lane $j$, whereas $p_n^{k}$ and $s_n^{k}$ represent the vehicles just behind and in front of vehicle $n$ in lane $k=j-1,j,j+1$, respectively.}\label{strada}
\end{figure}
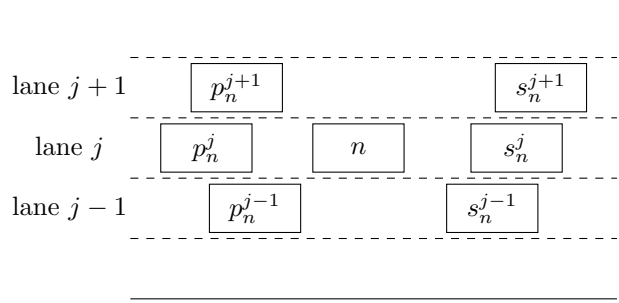

In the case of multi-lane traffic, vehicles travel along multiple lanes of length $L>0$ with the possibility to change lanes. Consider $J\in\mathbb{N}$ lanes and  denote by $j \in {1,\dots,J}$ the lane index. Given a vehicle $n$, traveling in the lane $j$, we denote with $p_n^{j}$ and $s_n^{j}$ the following and leading vehicles, respectively, in lane $j$. Further, $p_n^{k}$ and $s_n^{k}$, $k=j-1,j+1$, denote the vehicles in the two adjacent lanes; see Fig.~\ref{strada}. If    on lane $j'\in\{j-1,j+1\}$ there is a vehicle in the same   position of vehicle $n$ we choose to name it as $p_n^{j'}$. In the following the label  assigned to a  vehicle on the road   remains unchanged for  all time, even after a lane change occurs. 

Let $\Delta x^j_n = x_{s^j_n}-x_n$ be the headway between the $n$-th vehicle and its leading vehicle in lane $j$. The first order multi-lane microscopic model is written for $j\in \{1,\dots,J\}$ as a system of ordinary differential equations given by
\begin{equation}\label{bftl_2lane}
\begin{cases}
\dot x_n =V_j(\Delta x^j_n), \qquad n\in I_j(t), \\
\text{ + lane changing conditions}
\end{cases} \qquad \text{for } j=1,\dots,J,
\end{equation}
where $I_j(t)$ denotes the set of vehicles in lane $j$ at time $t$, ordered by the position. The function $V_j(h)$ prescribes the desired velocity  depending on the headway $h$. We assume $V_j$ to be a monotonic increasing function with $V_j(0)=0$ and  bounded by a lane-dependent maximum  value $V^{\max}_j$.

Equation~\eqref{bftl_2lane} is coupled to a set of  discrete rules that describe the occurrence of a lane change. We assume to deal with instantaneous lane changing events occurring from the current lane to one of the two adjacent lanes, the target lane. A vehicle performs a lane change based on two criteria~\cite{PiuPuppo2022,GongPiccoliVisconti}: an \emph{incentive} criterion, i.e.~a lane change may occur if local speed in the target lane is greater than the local speed in the current lane, and a \emph{safety} criterion, i.e.~a vehicle moves to a target lane if it is possible to keep a safe distance $d_s>0$ in order to avoid collisions. In this way the lane changing conditions from lane $j$ to a lane $j'\in\{j-1,j+1\}$ are written as
\begin{equation}\label{microLCC}
\begin{cases}
V_{j'}(x_{s^{j'}_n}-x_n)>V_j(x_{s^{j}_n}-x_n) & \text{(incentive criterion)} \\[1ex]
x_{s^{j'}_n}-x_n>l+d_{s} \quad \text{and} \quad x_n-x_{p_n^{j'}}>l+d_{s} & \text{(safety criterion)} \end{cases}
\end{equation}
Note that if the current lane is either $j=1$ or $j=J$, the target lane is one, namely lane $2$ or lane $J-1$, respectively. Instead, when the current lane is $j\in\{2,\dots,J-1\}$ there are two target lanes $j-1$ and $j+1$. If both satisfy the lane changing rules we assume that a vehicle moves to the most advantageous one, in terms of the velocity. If both changes are possible and the vehicle can gain the same speed in both lanes we choose to prioritize the change to the left. Moreover, after a lane change a vehicle maintains the same position as in its original lane (only the index of the lane is updated).

For the purposes of this paper, of particular interest is the discussion of the steady states of~\eqref{bftl_2lane}.

\begin{definition} \label{def:ss:micro}
We say that system~\eqref{bftl_2lane}-\eqref{microLCC} is \emph{locally at equilibrium} if in each lane vehicles are equally spaced and, thus, move with  constant velocity. Namely, for any fixed $j\in\{1,\dots,J\}$ we have
$
    \Delta x_n^j = C_j \, \forall n\in I_j,
$ where $C_j > 0$ are lane-dependent constants. Moreover, we say that system~\eqref{bftl_2lane}-\eqref{microLCC} is \emph{globally at equilibrium} if all lanes are locally at equilibrium and there is no lane changing, i.e.,~\eqref{microLCC} is violated for any fixed $j\in\{1,\dots,J\}$ and $j'\in\{j-1,j+1\}$.
\end{definition}

\begin{remark}
    For the purposes of this paper, without loss of generality in the following we consider a multi-lane ring road of length $L > 0$. Thus, system~\eqref{bftl_2lane} is endowed with periodic boundary conditions on the domain $\left[0,L\right]$ $\forall j=1,\dots,J$.
\end{remark}

\subsection{Evolution of the local density: a discrete description}
As a first step to derive a macroscopic model including the microscopic rules~\eqref{bftl_2lane}-\eqref{microLCC} we provide a definition for the local density. We assume that on each lane the maximum value of the density, $\rho^{\max}$, is reached when the vehicles are equally spaced with a bumper-to-bumper distance equal to $d_s$. Imposing $\rho^{\max}=1$, we define the local density as follows.

\begin{definition}[local density] The local density of the $n$-th vehicle in lane $j$ at time $t$ is given by:
\begin{equation}\label{definizione_densita}
	\rho^{(n)}_j(t):=\frac{l+d_s}{\Delta x^j_n}.
\end{equation}
\end{definition}

\begin{remark}
With the previous definition the density is dimensionless. 
\end{remark}

In the following, we analyze the behavior of the local density after a single lane change. The vehicle $n$ moves from the current lane $j$ to the target lane $j'\in\{j-1,j+1\}$ assuming that conditions~\eqref{microLCC} are satisfied at  time $t$. At time $t^-$, before the lane change occurs, the local density around the vehicle $p_n^{j'}$ is given by
\begin{equation*}
	\rho_{j'}^{(p^{j'}_n)} (t^-)= \frac{l+d_s}{x_{s^{j'}_n}-x_{p^{j'}_n}},
\end{equation*} 
whereas  after the lane change has occurred at time $t^+$, the local density in front of the vehicle $p_n^{j'}$ becomes
\begin{equation*}
\rho_{j'}^{(p^{j'}_n)} (t^+)= \frac{l+d_s}{\tilde{x}-x_{p^{j'}_n}},
\end{equation*}
where $\tilde{x}$ denotes the new position of vehicle $n$ in lane $j'$. We write $\tilde{x}$ as a convex combination of the positions of the follower and the leading vehicles in the new lane $j'$. Recalling the safety criterion in~\eqref{microLCC}, we have
\begin{equation*}\label{combinazione}
	\tilde{x}=\lambda (x_{s^{j'}_n}-d_s-l) + (1-\lambda) (x_{p^{j'}_n}+l+d_s), \qquad \lambda\in[0,1].
\end{equation*}
The parameter $\lambda$ can be modeled as a constant or as a function of the local density of vehicle $n$ in its current lane $j$. This allows us to describe the following scenario: if the density $\rho_{j}^{(n)}(t^-)$ is large, it is more likely that the position $\tilde{x}$ of vehicle $n$ in the new lane $j'$ is closer to the follower vehicle $p_n^{j'}$, because its velocity will be relatively small:
\begin{equation*}\label{gammadefini}
\lambda=\lambda(\rho_{j}^{(n)}(t^-)) \quad \text{s.t.} \quad \lambda:[0,\rho^{\text{max}}]\to[0,1]\quad \text{and}  \quad \lambda'(\cdot)<0, \,\, \lambda(0) =1,\,\,\lambda(\rho^{\text{max}}) =0.
\end{equation*}
The simplest choice is the linear function $ \lambda(\rho)=1-\frac{\rho}{\rho^{\text{max}}}. $

The density increment on lane $j'$ after the lane change of vehicle $n$ from lane $j$ can be computed as the difference:
\begin{align*} 
\rho_{j'}^{(p^{j'}_n)}(t^+)-\rho_{j'}^{(p^{j'}_n)}(t^-)
&= \left(\frac{1}{\lambda(\rho_j^{(n)}(t^-)) + (1-2\lambda(\rho_j^{(n)}(t^-)))\rho_{j'}^{(p^{j'}_n)}(t^-)}-1 \right)\rho_{j'}^{(p^{j'}_n)}(t^-) \\
&=:  A\left(\rho_j^{(n)}(t^-),\rho_{j'}^{(p^{j'}_n)}(t^-)\right)\rho_{j'}^{(p^{j'}_n)}(t^-),
\end{align*}
where the function $A:(0,\rho^{\text{max}}] \times [0,\mu) \to (0,+\infty)$ expresses the amplification factor of the density in the target lane. Here, $\mu<\rho^{\text{max}}$ is a critical density to be detailed later. As we will see $\mu$ is the maximum density in the target lane allowing lane change.  

Summarizing, the new local density in lane $j'$ when a vehicle $n$ moves from lane $j$ to lane $j'$ is
\begin{equation*}
\rho_{j'}^{(p^{j'}_n)}(t^+)=
\rho_{j'}^{(p^{j'}_n)}(t^-)+ A\left(\rho_{j}^{(n)}(t^-),\rho_{j'}^{(p^{j'}_n)}(t^-)\right)\rho_{j'}^{(p^{j'}_n)}(t^-).
\end{equation*}

\begin{remark} In the case of zero densities in the target lane a different source term must be included. In such case, after a lane change, the density on the target lane increases from 0 to 
\begin{equation*}
\rho_{j'}^{(n)}(t^+)= \frac{l+d_s}{L}>0.
\end{equation*} 
\end{remark}

\section{Derivation of the first order multi-lane macroscopic model} \label{sec:derivation}
We present the derivation of the first order  multi-lane macroscopic model obtained as a continuous limit of the first order microscopic model with lane changing conditions introduced in the previous section.

\subsection{Evolution of the local density: a continuous  description} \label{sec:macro:evolution}
The previous description of the local discrete density allows us to derive a macroscopic scale formulation for the evolution of the density in each lane. Using a piecewise constant reconstruction of the local densities we define the macroscopic density as follows.

\begin{definition}[macroscopic density] The macroscopic density is given by
\begin{equation}\label{piece}
\rho_{j}(x,t):=\rho_{j}^{(n)}(t), \quad n \in I_j, \quad x \in [x_n,x_{s^{j}_n}) \subset \left[0,L\right],
\end{equation}
where $I_j=I_j(t)$ is the label set of the vehicles in lane $j$ and $L > 0$ is the length of the road. 
\end{definition}

Thus the macroscopic density is a piecewise constant function interpolating  the microscopic data. 
Since mass is conserved, the space-time evolution of $\rho_j=\rho_j(x,t)$ without lane changing must satisfy, weakly, a scalar conservation law of the form 
\begin{equation}\label{lwr}
\de_t \rho_j + \de_x f_j(\rho_j) = 0, \quad \text{on } (x,t)\in\left[0,L\right]\times \mathbb{R}^+, \quad \forall j\in \{1,\dots,J\},
\end{equation}
where the function $f_j(\cdot)$ is the flux function of lane $j$. For first order models
$$
    f_j(\rho_j)=\rho_j v_j(\rho_j).
$$
Here, $v_j$ is the macroscopic velocity  in lane $j$,  and it will be defined by the desired velocity function of the microscopic model~\eqref{bftl_2lane} as
\begin{equation} \label{eq:micro:macro:speed}
    v_j(\rho):= V_j \left(\frac{l+d_s}{\rho}\right).
\end{equation}

Let us discuss in detail the evolution of the densities $\rho_j$ in the presence of lane changing.

From~\eqref{lwr}, using  a Taylor expansion for the time derivative with $\Delta t>0$ we obtain  
\begin{equation*}
    \rho_j(x,t+\Delta t) = \rho_j - \Delta t \de_x f_j(\rho_j)+ o(\Delta t) \qquad \mbox{ for } j=1,\dots,J.
\end{equation*}
In order to study the effects of the lane changes on the macroscopic evolution of the density $\rho_j$, we concentrate on considering  the lane changes from lane $j'$ to lane $j$. Observe that $j'$ can be either $j-1$ or $j+1$. Let $X\sim\text{Ber}(\pi)$ be a random variable with Bernoulli distribution, parametrized by $\pi$. Let $X=1$ denote the event ``the vehicle jumps to a new lane in a unit of time'', while for $X=0$ there is no lane change in the unit of time. Then $\mathbb{P}(X=1)=\pi $, whereas $\mathbb{P}(X=0)=1-\pi$. 
In this setting, the density changes stochastically depending on the outcome of $X$ as: 
\begin{equation}\label{casi}
\tilde{\rho}_j(x,t+\Delta t;X)= 
\rho_j + \Delta t \left(-\de_x f_j(\rho^{n}_j) + \ X \nu A(\rho_{j'},\rho_{j})\rho_j\right)+ o(\Delta t),
\end{equation}
where we have introduced the term $\nu>0$ that models the frequency of lane changes, with dimension one over time. For the purpose of this study we assume that $\nu$ is constant and it is equal for all lanes.

We  compute the time evolution of the density as  expected value of \eqref{casi}
\begin{equation*} 
\begin{split}
\rho_j (x,t+\Delta t)
& =\mathbb{E}[\tilde{\rho}_j(x,t+\Delta t;X)] = \rho_j + \Delta t \left(- \de_x f_j(\rho_j) + \pi \nu A(\rho_{j'},\rho_{j})\rho_j \right)+ o(\Delta t),
\end{split}
\end{equation*}
and finally,
\begin{equation*}
\frac{\rho_j (x,t+\Delta t)- \rho_j(x,t)}{\Delta t}= - \de_x f_j(\rho_j(x,t)) + {\pi}{\nu} A(\rho_{j'}(x,t),\rho_j(x,t))\rho_j(x,t)+o(1).
\end{equation*}
Renaming the probability $\pi$ as $\pi^{j' \to j}$ to emphasize the lane change from lane $j'$ to lane $j$, we obtain for $\Delta t \to 0$ the continuous description
\begin{equation} \label{eq:macro:jpTOj}
\de_t \rho_{j} + \de_x f_j(\rho_{j}) = {\pi^{j' \to j}}{\nu} A(\rho_{j'},\rho_{j})\rho_{j}, \quad \text{on } (x,t)\in\left[0,L\right]\times \mathbb{R}^+, \quad \forall j=1,\dots,J,
\end{equation}
with
\begin{equation*} \label{eq:A}
A(\rho_{j'},\rho_{j})= \left(\left(\lambda(\rho_{j'}) + (1-2\lambda(\rho_{j'}))\rho_{j}\right)^{-1}-1 \right),
\end{equation*}
where the source term describes the increase of the density in lane $j$ due to lane changing maneuvers from lane $j'$ to lane $j$.

\begin{remark}
We remark that in our micro-macro derivation in the case of a constant $\lambda$ the source term of the evolution equation of the density in lane $j$ is directly proportional to $\rho_{j}$ itself, whereas in other models, e.g.,~\cite{HoldenRisebro2019,SongKarni2019}, the source term is typically directly proportional to the density of the adjacent lanes.
\end{remark}

\subsection{Lane changing conditions in the macroscopic scale} \label{sec:macro:lcc}
In order to complete the macroscopic scale description of the evolution of the density in lane $j$ due to maneuvers from lane $j'$ to lane $j$, we  relate the lane changing probability with the microscopic lane changing conditions~\eqref{microLCC}. We start recalling the incentive and safety criteria, rewriting them in terms of the previously introduced macroscopic variables.

The incentive condition for a lane changing from $j'$ to $j$, i.e.
\begin{equation*}
	V_{j}(x_{s_n^{j}}-x_n) > V_{j'}(x_{s_n^{j'}}-x_n),
\end{equation*}
is given in the macroscopic setting using~\eqref{eq:micro:macro:speed} as
\begin{equation*}
	v_{j}\left(\frac{l+d_s}{x_{s_n^{j}}-x_n}\right) >v_{j'}\left(\frac{l+d_s}{x_{s_n^{j'}}-x_n}\right),
\end{equation*}
which implies 
\begin{equation*}
	v_{j}(\rho^{(n)}_{j}) > v_{j'}(\rho^{(n)}_{j'}).
\end{equation*}
Notice that $\rho^{(n)}_{j}$ is an \emph{expected} density since it describes the local density in front of vehicle $n$ in case it moves to lane $j$. In order to rewrite the incentive condition at the macroscopic level, instead of considering this expected density we take into account the current density in the target lane. This corresponds to the density observed by the driver before performing a lane changing. Indeed, matter-of-factly, drivers consider moving to another lane by evaluating whether the speed in that lane is higher than the speed in their current lane. 

Hence  the following condition on the macroscopic speeds is given:
\begin{equation}\label{cond1}
	v_{j}(\rho_{j}) > v_{j'}(\rho_{j'}).
\end{equation}
Therefore, no lane change occurs if the flow speed in the target lane $j$ is lower than the flow speed in the current lane $j'$.
On the other hand, the safety criterion for a lane change from $j'$ to $j$ can be obtained summing the two microscopic inequalities in~\eqref{microLCC}
\begin{equation*}
x_{s_n^{j}}-x_{p_n^{j}}>2(l+d_s),
\end{equation*}
and since $x_{s_n^{j}}-x_{p_n^{j}}>0$, this gives
\begin{equation} \label{eq:mu:def}
\frac{l+d_s}{x_{s_n^{j}}-x_{p_n^{j}}}<\frac{1}{2}=:\mu \text{   (critical density)}. 
\end{equation}
Thus the safety criterion reduces to 
\begin{equation*} 
\rho^{(p_n^{j})}_{j} < \mu,
\end{equation*}
and in terms of the macroscopic density, we obtain
\begin{equation}\label{cond2}
\rho_{j} < \mu.
\end{equation}
%
Therefore, no lane change occurs if $\rho_{j}\geqslant \mu$, and we expect that the probability of changing lanes increases as $\mu -\rho_{j}$ becomes larger.

Summarizing the previous discussion, a lane change can occur only if conditions \eqref{cond1}-\eqref{cond2} are satisfied. Therefore, the probability of a lane change from lane $j'$ to lane $j$, denoted by $\pi^{j' \to j}$, in a macroscopic setting, is derived from the lane changing conditions as follows.

Let $\mathbf{1}_{lc}(\rho_{j'},\rho_j)$ be the indicator function from lane changing from lane  $j'$ to lane $j$. Then
\begin{equation*}
\begin{split}
\mathbf{1}_{lc}(\rho_{j'},\rho_j) & =\begin{cases}
1 & \text{   if conditions  } \begin{cases}
v_{j}(\rho_{j}) > v_{j'}(\rho_{j'}) \\
\rho_{j} < \mu
\end{cases} \text{   are met,} \\
 0 & \text{   otherwise   }
\end{cases}\\
& = \max\left\{0, \min\left\{ \sign\left(v_{j}(\rho_{j}) - v_{j'}(\rho_{j'})\right), \sign\left(\mu -\rho_{j}\right) \right\}\right\}.
\end{split}
\end{equation*}
Therefore the probability $\pi^{j' \to j}$ is defined as
\begin{equation*}
\begin{aligned}
\pi^{j' \to j}(\rho_{j'},\rho_{j})= g(\rho_{j})\, \mathbf{1}_{lc}(\rho_{j'},\rho_j),
\end{aligned}
\end{equation*}
where the function $g:[0,\mu]\to[0,1]$ models the percentage of mass going from lane $j'$ to lane $j$.
Here, $g(\cdot)$ is a monotonically decreasing function of its argument. An example is provided by the linear function $$ g(\rho)=1-2\frac{\rho}{\rho^{\text{max}}}. $$

\begin{remark}
We point out that so far we have focused on the evolution of the density in lane $j$ when a lane change from $j'\in\{j-1,j+1\}$ to $j$ occurs. This causes an increase of $\rho_j$ and corresponds to a \emph{gain term}. Similarly, there is a decrease of $\rho_{j}$ which is modelled by the same source term in~\eqref{eq:macro:jpTOj} with a negative sign,  giving a \emph{loss term}. Hence we conserve the mass of all lanes.\end{remark}

\subsection{Final form of the model}
Here, we summarize the structure of the complete model. The derivation in Subsection~\ref{sec:macro:evolution} and Subsection~\ref{sec:macro:lcc} leads to the following final first order macroscopic model for a multi-lane road with $J$ lanes:
\begin{equation}\label{modello}
	\de_t \rho_{j} + \de_x f_j(\rho_{j}) = S(\rho_j,\rho_{j'}), \quad \text{on } (x,t)\in\left[0,L\right]\times \mathbb{R}^+, \quad \forall j=1,\dots,J,
\end{equation}
where the source term is defined as
\begin{equation*} \label{eq:source}
S(\rho_j,\rho_{j'}) = \sum\limits_{j' \in T_j}\pi^{j' \to j}(\rho_j,\rho_{j'}) \nu A(\rho_{j'},\rho_{j})\rho_{j} - \sum\limits_{j' \in T_j} {\pi^{j \to j'}}(\rho_{j'},\rho_{j}) \nu A(\rho_{j},\rho_{j'})\rho_{j'},
\end{equation*}
with
\begin{equation*} \label{eq:gamma:and:A}
\begin{aligned}
\pi^{h \to k}(\rho_h,\rho_k)  
&= g(\rho_k)\, \mathbf{1}_{lc}(\rho_{h},\rho_k), \\
\mathbf{1}_{lc}(\rho_{h},\rho_k) &= \max\left\{0, \min\left\{ \sign\left(v_{k}(\rho_{k}) - v_{h}(\rho_{h})\right), \sign\left(\mu -\rho_{k}\right) \right\}\right\}, \\
A(\rho_{h},\rho_{k}) &= \left(\left(\lambda(\rho_{h}) + (1-2\lambda(\rho_{h}))\rho_{k}\right)^{-1}-1 \right),
\end{aligned}
\end{equation*} 
 $ \forall k\in\{1,\dots,J\}$ and $h\in T_j$. The set $T_j$ of the target lanes of a given lane $j$ is
\begin{equation*}
T_j=
\begin{cases}
\{2\} & \mbox{ if }j=1,\\
\{j-1,j+1\} & \mbox{ if } j=2,\dots,J-1,\\
\{J-1\} & \mbox{ if }j=J.
\end{cases}
\end{equation*}
Furthermore $\lambda(\rho)=1-\frac{\rho}{\rho^{\text{max}}}, \quad \mbox{and} \quad g(\rho)=1-2\frac{\rho}{\rho^{\text{max}}}. $
System~\eqref{modello} is endowed with given initial conditions $\rho_j^0(x):=\rho_j(x,t=0)$, and suitable boundary conditions.

\section{Steady state analysis} \label{sec:analysis}

This section is devoted to the discussion of the stability of the steady states of the macroscopic multi-lane model~\eqref{modello}. We are interested in the particular class of steady states defined as follows.

\begin{definition} \label{def:ss}
A steady state $\{\overline{\rho}_j(x)\}_{j=1}^J$ of the macroscopic multi-lane model \eqref{modello} is characterized by
$$
  \overline{\rho}_j(x) = \overline{\rho}_j\in (0,\rho^{\max}] \ \text{ and } \ S(\overline{\rho}_j,\overline{\rho}_{j'}) = 0 \quad \forall j=1,\dots,J, \quad j'\in T_j,
$$
i.e.~in each lane the density is constant and no lane changes occur.
\end{definition}

This corresponds to the global equilibrium of the microscopic model defined in Definition~\ref{def:ss:micro}.

\subsection{Explicit characterization of the equilibria} \label{ssec:charact:eq}
For the sake of simplicity, we consider the model \eqref{modello} in the case of a two-lane road, i.e.~$J=2$, described by the system of two balance laws:
\begin{equation}\label{duecorsie}
\begin{cases}
\de_t \rho_{1} + \de_x f_1(\rho_{1}) = {\pi^{2 \to 1}}{\nu} A(\rho_{2},\rho_{1})\rho_{1} - {\pi^{1 \to 2}}{\nu}A(\rho_{1},\rho_{2})\rho_{2}, \\

\de_t \rho_{2} + \de_x f_2(\rho_{2}) = {\pi^{1 \to 2}}{\nu} A(\rho_{1},\rho_{2})\rho_{2}- {\pi^{2 \to 1} }{\nu}A(\rho_{2},\rho_{1})\rho_{1}. \\
\end{cases}
\end{equation}

Since $\nu$ is a scaling factor, we assume here $\nu=1$ and omit to write it in the following analysis.

\begin{figure}[t]
\centering
\begin{tikzpicture}[scale=0.8 ]
\draw[very thick,name path=v2,blue] (0,4)--(4,0);
\draw[very thick,name path=v1,red] (4,0)--(0,3);
\draw[->] (0,0) -- (5,0)node[below]{$\rho$}; 
\draw[->] (0,0) -- (0,5)node[left]{$v$};
\draw (2,0)node[below]{$\mu$};
\draw (4,0)node[below]{$\rho^{\text{max}}$};
\draw (0,4)node[left]{$v_2^{\text{max}}$};
\draw (0,3)node[left]{$v_1^{\text{max}}$};
\draw [dashed] (2,0)--(2,2);
\draw [dashed] (2,2)--(0,2);
\draw (0,2)node[left]{$v_2^{\mu}$};
\draw [dashed] (2,-3/2+3)--(0,-3/2+3);
\draw (0,-3/4*2+3)node[left]{$v_1^{\mu}$};
\draw[dotted] (4/3,2)--(4/3,0);
\draw (4/3,0)node[below]{$\rho_1^{\mu}$};
\draw[dotted] (1+3/2,-3/2+3)--(1+3/2,0);
\draw (1+3/2,0)node[below]{$\rho_2^{\mu}$};
\draw [dashed] (2,-3/2+3)--(1+3/2,-3/2+3);
\end{tikzpicture}
\caption[An example with linear velocity function.]{An example with linear velocity function. Blue: $v_2$; red: $v_1$.}\label{img:notazioni}
\end{figure}
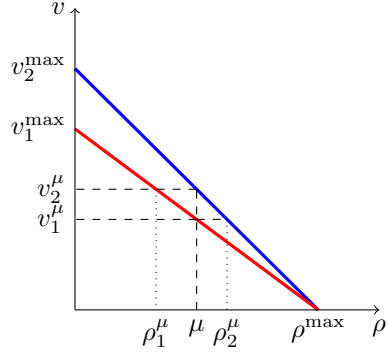

We denote $\overline{\rho}_1$ and $\overline{\rho}_2$ the steady state in lane 1 and lane 2, respectively. Furthermore, we assume that the two speed functions $v_1(\rho),v_2(\rho)$ are monotonically decreasing with respect to the density $\rho$, and are chosen as
\begin{equation} \label{eq:speed:functions}
v_1(\rho)<v_2(\rho),  \quad \rho\in[0,\rho^{\max}) \quad \text{ and } \quad \begin{aligned}
v_1(\rho^{\max})&=v_2(\rho^{\max})=0, \\ v_j(0)&=v_j^{\max}, \quad j=1,2.
\end{aligned}
\end{equation}
In the following, we denote  by $v_j^\mu=v_j(\mu)$ the flow speed on the $j$-th lane corresponding to the critical density $\mu$ introduced in~\eqref{eq:mu:def}. In addition, we use the notation
$$
\rho_j^\mu:=(v_j)^{-1}(v_{j'}^\mu), \quad j,j'=1,2, \ j\neq j',
$$
to denote the unique value of the density in lane $j$ whose corresponding flow speed is equal to the speed of lane $j'$ at the critical density. See Fig.~\ref{img:notazioni} for an example with the linear velocities 
\begin{equation} \label{lineari}
\begin{aligned}
v_j(\rho) &= v_j^{\max}\left(1-\frac{\rho}{\rho^{\max}}\right), \quad j=1,2, \\
0&<v_1^{\max} < v_2^{\max}.
\end{aligned}
\end{equation}

We recall that lane changes are possible when both the incentive and the safety criteria are satisfied. Depending on the realization of the two criteria, the following four characterizations of the steady states given in Definition~\ref{def:ss}, where no lane changes occur, arise.

\begin{itemize}
\item[(A)] The density in both lanes is smaller than the critical density and the lane speeds are equal, namely
\begin{equation} \label{eq:eq:A}
v^{\text{eq}}:=v_1(\overline{\rho}_1)=v_2(\overline{\rho}_2), \qquad \text{ and }\qquad  0\leqslant \overline{\rho}_1 < \mu \,\,\land  \,\, 0\leqslant \overline{\rho}_2< \mu.
\end{equation}
\item[(B)] At least one of the two lanes has density greater than the critical value and the lane speeds are equal, namely
\begin{equation} \label{eq:eq:B}
v^{\text{eq}}:=v_1(\overline{\rho}_1)=v_2(\overline{\rho}_2), \qquad \text{ and }\qquad  \overline{\rho}_1\geqslant \mu \,\,\lor  \,\,  \overline{\rho}_2\geqslant \mu.
\end{equation}
\item[(C)] The density in both lanes is greater than the critical value and the lane speeds are different, namely
\begin{equation} \label{eq:eq:C}
v_1(\overline{\rho}_1)\neq v_2(\overline{\rho}_2), \qquad \text{ and }\qquad  \overline{\rho}_1 \geqslant \mu \,\,\land  \,\,  \overline{\rho}_2 \geqslant \mu.
\end{equation}
\item[(D)] The density in lane 1 is smaller than the critical value, whereas the density in lane 2 is greater, and the flow in lane 1 is slower than the flow in lane 2, namely
\begin{equation} \label{eq:eq:D}
v_1(\overline{\rho}_1)< v_2(\overline{\rho}_2), \qquad \text{ and }\qquad  \overline{\rho}_1< \mu \,\,\land  \,\,  \overline{\rho}_2 \geqslant \mu.
\end{equation}
\end{itemize}

To these equilibria we also add the case in which all the mass is in the fastest lane and lane changes do not occur because the incentive criterion is not satisfied.
\begin{itemize}
\item[(E)] The density in lane 1 is zero, whereas the density in lane 2 is less than or equal to $v_2^{-1}( v_1^{\text{max}})$, namely
\begin{equation} \label{eq:eq:E}
  \overline{\rho}_1 =0 \,\,\land  \,\,  \overline{\rho}_2\leqslant v_2^{-1}( v_1^{\text{max}}).
\end{equation}
\end{itemize}

\subsection{Linear stability analysis}

In the following, we assume that system~\eqref{duecorsie} is at equilibrium, according to Definition~\ref{def:ss}, and we consider a uniform perturbation in space on the two lanes such that the total mass is conserved. More precisely, the perturbed state $(\rho_1,\rho_2)$ of the equilibrium $(\overline{\rho}_1,\overline{\rho}_2)$ is given by
\begin{equation}\label{perturbazione_modello}
\begin{cases}
\rho_1=\overline{\rho}_1+\epi \rho\\
\rho_2=\overline{\rho}_2-\epi \rho
\end{cases} \quad \text{with } \ \rho=\rho(t), \ \rho(0)=\rho_0>0
\end{equation}
with $|\epi|\ll 1$. The aim is to investigate the stability of the five types of steady states previously described. In the following we divide the equilibria class (B) into (B1) with $0\leqslant v^{\text{eq}}\leqslant v_1^\mu$ and  (B2) with $v_1^\mu < v^{\text{eq}}\leqslant v_2^\mu$. The equilibria of type (B1) can be seen as a particular case of equilibria (C) with equal velocities $v_1(\overline{\rho}_1)= v_2(\overline{\rho}_2)$. Therefore, for simplicity of exposition, they will  now be  included in the equilibrium class (C).

The left column of Table~\ref{tabella_convergenze} reports all the possible types of steady states of the two-lane model~\eqref{duecorsie}, which are characterized in the equations from~\eqref{eq:eq:A} to~\eqref{eq:eq:E}. When these are perturbed according to~\eqref{perturbazione_modello}, system~\eqref{duecorsie} converges, for large times, to the corresponding steady states listed in the right column of Table~\ref{tabella_convergenze}. In particular, we list the cases in which the perturbed system is driven towards the initial steady state, and the cases in which the system converges to stationary states, which are different from the initial one.

\begin{table}[t]
\centering
\caption{Left: equilibrium states of system~\eqref{duecorsie} as defined in~\eqref{eq:eq:A}-\eqref{eq:eq:E} and \eqref{defb2}. Right: equilibrium states of system~\eqref{duecorsie} when the corresponding states on the left column are perturbed according to~\eqref{perturbazione_modello}.}\label{tabella_convergenze}
\begin{tabular}{ |c|c|  }
\hline
Initial equilibrium & Final equilibrium \\
\hline \hline
(A) & (A) \\
(B2) or (C) or (D) & (B2) or (C) or (D) \\
(E) & (E) \\
\hline
\end{tabular}
\end{table}

The following definition allows us to characterize the steady states of system~\eqref{duecorsie}.

\begin{definition}[taxonomy of the equilibrium states] \label{def:stability}
    We say that an equilibrium $(\overline{\rho}_1,\overline{\rho}_2)$ of system~\eqref{duecorsie} according to Definition~\ref{def:ss} is
    \begin{itemize}
    \item \emph{asymptotically stable} with respect to the space homogeneous perturbation~\eqref{perturbazione_modello} if there exists $\epi\in\mathbb{R}$ such that $(\rho_1,\rho_2)\to(\overline{\rho}_1,\overline{\rho}_2)$ as $t\to\infty$;
    \item \emph{globally asymptotically stable} with respect to the space homogeneous perturbation~\eqref{perturbazione_modello} if $\forall\,\epi\in\mathbb{R}$ one has $(\rho_1,\rho_2)\to(\overline{\rho}_1,\overline{\rho}_2)$ as $t\to\infty$;
    \item \emph{marginally stable} with respect to the space homogeneous perturbation~\eqref{perturbazione_modello} if for a given $\epi\in\mathbb{R}$ one has $(\rho_1,\rho_2)\to(\tilde{\rho}_1,\tilde{\rho}_2)$ as $t\to\infty$, where $(\tilde{\rho}_1,\tilde{\rho}_2)$ is an equilibrium different from $(\overline{\rho}_1,\overline{\rho}_2)$.
    \end{itemize}
\end{definition}

By the definition of a marginally stable equilibrium we mean a trajectory starting from the perturbed state which approaches a new equilibrium. In other words, those perturbed states do not result in limit cycles.

Using Definition~\ref{def:stability} we state the nature of all steady states of system~\eqref{duecorsie} in the following theorem, whose proof is reported in the appendix.

\begin{theorem}[stability for the macroscopic two-lane model] \label{teorema}
Consider a space homogeneous perturbation as~\eqref{perturbazione_modello}. Then the steady states $(\overline{\rho}_1,\overline{\rho}_2)$ of the two-lane model~\eqref{duecorsie} can be characterized as follows according to Definition~\ref{def:stability}:
\begin{itemize}
    \item equilibrium solutions of type (A) and (E) are globally asymptotically stable;
    \item equilibrium solutions of type (B2) and (C) with $\overline{\rho}_1=\mu$ are asymptotically stable if $\epi<0$, and marginally stable otherwise;
    \item equilibrium solutions of type (C) and (D) with $\overline{\rho}_2=\mu$ are asymptotically stable if $\epi>0$, and marginally stable otherwise;
    \item equilibrium solutions of type (C) with $\overline{\rho}_1> \mu$ and $\overline{\rho}_2> \mu$  and type (D) with  $\overline{\rho}_2>\mu$ are marginally stable.
\end{itemize}
\end{theorem}

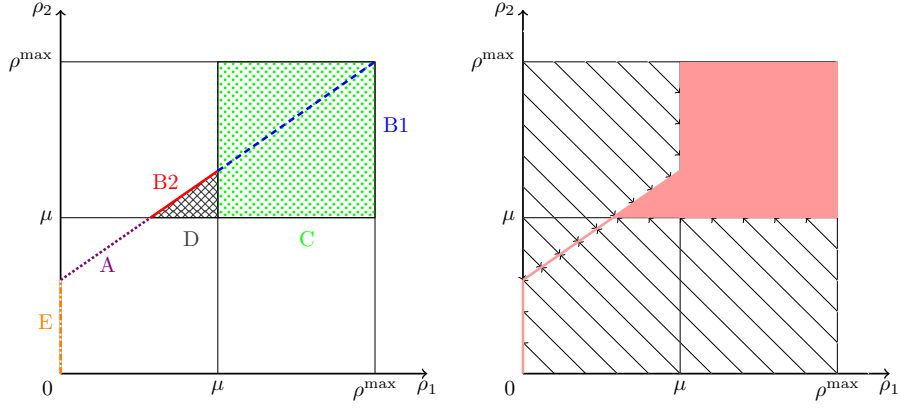
\begin{figure}[t]
\centering
\resizebox{\columnwidth}{!}{%
\begin{tikzpicture}
\begin{scope}[scale=0.85]
\draw[thick, ->](0,0)--(7,0);
\draw[thick, ->](0,0)--(0,7);
\path[below](7,0)node{$\rho_1$};
\path[left](0,7)node{$\rho_2$};
\path[below](6,0)node{$\rho^{\max}$};
\path[left](0,6)node{$\rho^{\max}$};
\path[below](3,0)node{$\mu$};
\path[left](0,3)node{$\mu$};
\path[anchor=north east](0,0)node{$0$};
\draw[name path= mu1, very thin] (0,3)--(6,3);
\draw[name path= mu2, very thin] (3,0)--(3,6);
\draw[very thin] (0,6)--(6,6);
\draw[very thin] (6,0)--(6,6);
\draw[white, name path=eq](0, 0.3*6)--(6,6); 
\path [name intersections={of=eq and mu1,by=E}];
\path [name intersections={of=eq and mu2,by=F}];
\filldraw[pattern=\patP,pattern color=darkgray] (E)--(F)--(3,3)--cycle;
\filldraw[pattern=\patUnoMenoPmezzi, pattern color=green] (3,3)--(3,6)--(6,6)--(6,3)--cycle;
\draw[very thick,violet,densely dotted] (0,0.3*6)--(E);
\draw[very thick,red] (E)--(F);
\draw[very thick,blue,densely dashed] (F)--(6,6);
\draw[thick,white] (0,0)--(0,0.3*6);
\draw[very thick,orange,densely dash dot] (0,0)--(0,0.3*6);
\path[left](0,1)node{\textcolor{orange}{E}};
\path[left](1.2,2.1)node{\textcolor{violet}{A}};
\path[left](2.4,3.7)node{\textcolor{red}{B2}};
\path[left](2.8,2.6)node{\textcolor{darkgray}{D}};
\path[left](5,2.6)node{\textcolor{green}{C}};
\path[left](6.8,4.8)node{\textcolor{blue}{B1}};
\end{scope}
\begin{scope}[xshift=7.5cm,scale=0.85]
\draw[thick, ->](0,0)--(7,0);
\draw[thick, ->](0,0)--(0,7);
\path[below](7,0)node{$\rho_1$};
\path[left](0,7)node{$\rho_2$};
\path[below](6,0)node{$\rho^{\max}$};
\path[left](0,6)node{$\rho^{\max}$};
\path[below](3,0)node{$\mu$};
\path[left](0,3)node{$\mu$};
\path[anchor=north east](0,0)node{$0$};
\draw[name path= mu1, very thin] (0,3)--(6,3);
\draw[name path= mu2, very thin] (3,0)--(3,6);
\draw[name path= a1, very thin] (0,6)--(6,6);
\draw[name path= a2, very thin] (6,0)--(6,6);
\draw[white, name path=eq](0, 0.3*6)--(6,6); 
\path [name intersections={of=eq and mu1,by=E}];
\path [name intersections={of=eq and mu2,by=F}];

\draw[->](6*0.2,0)--(0,6*0.2);
\draw[->](6*0.1,0)--(0,6*0.1);

\draw[white, name path=eq1](6*0.3,0)--(0,6*0.3);
\path [name intersections={of=eq1 and eq,by=E0}];
\draw[->](6*0.3,0)--(E0);
\draw[->](0,6*0.3)--(E0);

\draw[white, name path=eq1](6*0.4,0)--(0,6*0.4);
\path [name intersections={of=eq1 and eq,by=E1}];
\draw[->](6*0.4,0)--(E1);
\draw[->](0,6*0.4)--(E1);

\draw[white, name path=eq1](6*0.5,0)--(0,6*0.5);
\path [name intersections={of=eq1 and eq,by=E2}];
\draw[->](6*0.5,0)--(E2);
\draw[->](0,6*0.5)--(E2);

\draw[white, name path=eq1](6*0.6,0)--(0,6*0.6);
\path [name intersections={of=eq1 and eq,by=E3}];
\draw[->](6*0.6,0)--(E3);
\draw[->](0,6*0.6)--(E3);

\draw[white, name path=eq1](6*0.7,0)--(0,6*0.7);
\path [name intersections={of=eq1 and eq,by=E4}];
\draw[->](6*0.7,0)--(E4);
\draw[->](0,6*0.7)--(E4);

\draw[white, name path=eq1](6*0.8,0)--(0,6*0.8);
\path [name intersections={of=eq1 and eq,by=E51}];
\path [name intersections={of=eq1 and mu1,by=E52}];
\draw[->](6*0.8,0)--(E52);
\draw[->](0,6*0.8)--(E51);

\draw[white, name path=eq1](6*0.9,0)--(0,6*0.9);
\path [name intersections={of=eq1 and eq,by=E61}];
\path [name intersections={of=eq1 and mu1,by=E62}];
\draw[->](6*0.9,0)--(E62);
\draw[->](0,6*0.9)--(E61);

\draw[white, name path=eq1](6*1,0)--(0,6*1);
\path [name intersections={of=eq1 and eq,by=E71}];
\path [name intersections={of=eq1 and mu1,by=E72}];
\draw[->](6*1,0)--(E72);
\draw[->](0,6*1)--(E71);

\draw[white, name path=eq1](6*1.1,0)--(0,6*1.1);
\path [name intersections={of=eq1 and a1,by=A1}];
\path [name intersections={of=eq1 and a2,by=A2}];
\path [name intersections={of=eq1 and eq,by=E81}];
\path [name intersections={of=eq1 and mu1,by=E82}];
\draw[->](A1)--(E81);
\draw[->](A2)--(E82);

\draw[white, name path=eq1](6*1.2,0)--(0,6*1.2);
\path [name intersections={of=eq1 and a1,by=A3}];
\path [name intersections={of=eq1 and a2,by=A4}];
\path [name intersections={of=eq1 and mu2,by=E91}];
\path [name intersections={of=eq1 and mu1,by=E92}];
\draw[->](A3)--(E91);
\draw[->](A4)--(E92);

\draw[white, name path=eq1](6*1.3,0)--(0,6*1.3);
\path [name intersections={of=eq1 and a1,by=A5}];
\path [name intersections={of=eq1 and a2,by=A6}];
\path [name intersections={of=eq1 and mu2,by=E101}];
\path [name intersections={of=eq1 and mu1,by=E102}];
\draw[->](A5)--(E101);
\draw[->](A6)--(E102);

\draw[white, name path=eq1](6*1.4,0)--(0,6*1.4);
\path [name intersections={of=eq1 and a1,by=A7}];
\path [name intersections={of=eq1 and a2,by=A8}];
\path [name intersections={of=eq1 and mu2,by=E111}];
\path [name intersections={of=eq1 and mu1,by=E112}];
\draw[->](A7)--(E111);
\draw[->](A8)--(E112);

\fill[red!40!white] (E)--(F)--(3,3)--cycle;
\fill[red!40!white] (3,3)--(3,6)--(6,6)--(6,3)--cycle;
\draw[very thick,red!40!white] (0,0.3*6)--(E);
\draw[very thick,red!40!white] (E)--(F);
\draw[very thick,red!40!white] (F)--(6,6);

\draw[very thick,red!40!white] (0,0)--(0,0.3*6);
\end{scope}
\end{tikzpicture}}
\caption[Equilibria in the phase portrait.]{Left: plane $(\rho_1,\rho_2)$ showing the steady states of model~\eqref{duecorsie}. The violet dotted line represents the \textcolor{violet}{equilibria of type (A)}, the red solid line the \textcolor{red}{equilibria of type (B2)}, the blue dashed line the \textcolor{blue}{equilibria of type (B1)}, the orange dash-dotted line the \textcolor{orange}{equilibria of type (E)}, the region shaded with criss-crossing dark gray lines represents the \textcolor{darkgray}{equilibria of type (D)}, and the region shaded with green dots represents the \textcolor{green}{equilibria of type (C)}. Right: phase portrait $(\rho_1,\rho_2)$ of model~\eqref{duecorsie} where the light red zones represent the equilibria of the system.}\label{fasi}
\end{figure} 

It is possible to represent the set of the equilibrium solutions of the two-lane model~\eqref{bftl_2lane} in the plane $(\rho_1,\rho_2)$. 

To illustrate the behaviour of the different equilibria we consider the case in which the speeds are linear functions. Thus the equilibrium equation $v_1(\rho_1)=v_2(\rho_2)$ can be written as  
\begin{equation*}
\rho_2=\rho^{\max}\left( 1-\frac{v_1^{\max}}{v_2^{\max}} \left(1-\frac{\rho_1}{\rho^{\max}}\right) \right)
\end{equation*} that it is a straight line in the $(\rho_1,\rho_2)$ plane, clearly visible in the left panel of Fig.~\ref{fasi} by  the violet dotted equilibria (A), the red solid equilibria (B2), and the blue dashed segments equilibria (B1).

Moreover, we report  the phase portrait of system~\eqref{duecorsie} in the right panel of Fig.~\ref{fasi}. The trajectories are depicted with black lines and they converge toward the equilibrium solutions, which are identified by the light red zones. We observe that the trajectories are straight lines since the perturbation we are considering is homogeneous in space. 

For more general speed functions, the types (A) and (B) equilibria lie on a curve whose equation depends on the expressions of such speed functions.

\section{Numerical tests} \label{sec:numerics}

In this section we present four numerical tests in order to investigate the behavior of the multi-lane macroscopic model~\eqref{modello} using two and three lanes. In particular, the first test aims to show the consistency between the microscopic multi-lane model and its macroscopic limit presented in Section~\ref{sec:derivation}. The subsequent two tests are devoted to studying the effects of space homogeneous and local (i.e.~space non-homogeneous) perturbations of equilibrium states. Finally, in the last test we propose a practical application of the multi-lane macroscopic model which is used to study the flow of vehicles on a three-lane highway with a local lane closure, e.g.~due to a work zone. In all simulations we set $\nu=1$ s\textsuperscript{$-1$}.

In all simulations the macroscopic system~\eqref{modello} is numerically integrated with a first order finite volume scheme with Rusanov numerical flux, considering a uniform grid of size $\Delta x$ and adaptive time steps $\Delta t$ in order to satisfy the Courant-Friedrichs-Lewy stability condition. 

\subsection{Consistency between the micro and the macro models}

\begin{table}[ht]
	\centering
	\caption{Numerical experiments on the consistency of the macroscopic limit. Values of the initial conditions and final states of the microscopic and macroscopic model.\label{tab:consi}}
	\resizebox{\columnwidth}{!}{%
	\begin{tabular}{l|cc|cc}
		 & \multicolumn{2}{c|}{Lane 1} & \multicolumn{2}{c}{Lane 2} \\
		& $t=0$ & $t=T=100$ & $t=0$ & $t=T=100$ \\
		\hline \hline
		TEST 1 & & & &\\
		\hline\hline
		\multirow{2}{*}{Microscopic} & $N_1=150$ & $N_1=104$ & $N_2=30$ & $N_2=76$ \\
		& $\rho_1^{(n)}=1$ & $\rho_1^{(n)}=0.69$ & $\rho_2^{(n)}=0.2$ & $\rho_2^{(n)}=0.51$ \\
		\hline
		Macroscopic & $\rho_{1}(x,0)=1$ & $\rho_{1}(x,T)=0.70$ & $\rho_{2}(x,0)=0.2$ & $\rho_{2}(x,T)=0.50$ \\
		\hline\hline
		TEST 2& & & &\\		
		\hline\hline
		\multirow{2}{*}{Microscopic} & $N_1=100$ & $N_1=72$ & $N_2=50$ & $N_2=78$ \\
		& $\rho_1^{(n)}=0.67$ & $\rho_1^{(n)}=0.48$ & $\rho_2^{(n)}=0.33$ & $\rho_2^{(n)}=0.52$ \\
		\hline
		Macroscopic & $\rho_{1}(x,0)=0.66$ & $\rho_{1}(x,T)=0.50$ & $\rho_{2}(x,0)=0.33$ & $\rho_{2}(x,T)=0.50$ \\
		\hline
	\end{tabular}}
\end{table}

In order to investigate the consistency between the microscopic and the macroscopic models, we focus on a two-lane circular road of length $\tilde{L}=1500$~meters and with vehicles of length $\tilde{l}=5$~meters. In the numerical simulations we consider normalized parameters, namely $L=1$, $l=d_s=0.00\overline{3}$, and a final time $T=100$.

The initial conditions of the microscopic model~\eqref{bftl_2lane} are set in such a way that the lanes are in local   equilibrium as in Definition~\ref{def:ss:micro}, but lane changes are possible. Therefore, we initially fix a suitable number of vehicles in the two lanes and impose that vehicles are uniformly distributed on each lane. Precisely, let $N_j(0)$ be the number of vehicles in lane $j=1,2$ at time $t=0$, then \begin{equation}\label{ic_micro}
x_{s_n^j}^{(t_0=0)}-x_n^{(t_0=0)}=h_j=\frac{L}{N_j(0)} \quad \forall n\in I_j(0),
\end{equation} where $I_j(0)$ is the set of vehicles, ordered by the position, in lane $j$ at time $t=0$.

For the macroscopic model~\eqref{duecorsie}, the initial condition is computed using the corresponding initial local discrete densities defined in~\eqref{piece}. Thus
\begin{equation}\label{ic_macro}
\rho_j(x,0)=\sum_{n\in I_j(0)}\rho_j^{(n)}(t=0) \chi_{[x_n,x_{s_n^j})}(x), \,\, j=1,2,
\end{equation}
where $\chi_I(x)$ is the characteristic function of the set $I\subset \mathbb{R}$. We employ the linear density-speed relations given in~\eqref{lineari} as macroscopic speed functions. In particular, we choose $v_2^{\max}=1$ and $v_1^{\max}=0.7$. Furthermore $\Delta x =0.00\overline{3}$.

The microscopic system is solved with a fifth-order Runge-Kutta method. The results of the two tests are summarized in Table~\ref{tab:consi}. Here, we report the chosen initial conditions, namely the number of vehicles and the local densities for the microscopic model and the macroscopic density for the continuum model, as well as the final states obtained via numerical evolution of both models. Both simulations show the excellent agreement between the microscopic multi-lane model and its macroscopic limit, at final time $t=T$.

\subsection{Global  perturbation in space  of equilibrium solutions}
In this test we show  numerically the transition toward equilibrium of the model \eqref{duecorsie}, confirming the theoretical study on the linear stability of the steady states. Since the analytical results provide a detailed description of the behavior of perturbations of equilibrium solutions, we present here only two numerical examples. We consider an initial datum for the two-lane model~\eqref{duecorsie} given by $(\rho_1,\rho_2)=(\overline{\rho}_1+\epi\rho_0,\overline{\rho}_2-\epi\rho_0)$, where $(\overline{\rho}_1,\overline{\rho}_2)=(0.27,0.49)$. Considering again the macroscopic density-speed relation as in~\eqref{lineari}, we have that $(\overline{\rho}_1,\overline{\rho}_2)$ belongs to the set of steady states of type (A). Note that $(\rho_1,\rho_2)$ is a uniform space perturbation of the steady state, as given in~\eqref{perturbazione_modello}. In particular, we choose $\epi\rho_0=0.485$ in the first test and $\epi\rho_0=-0.27$ in the second test.

\begin{figure}[ht]
\centering
\includegraphics[width=0.49\textwidth]{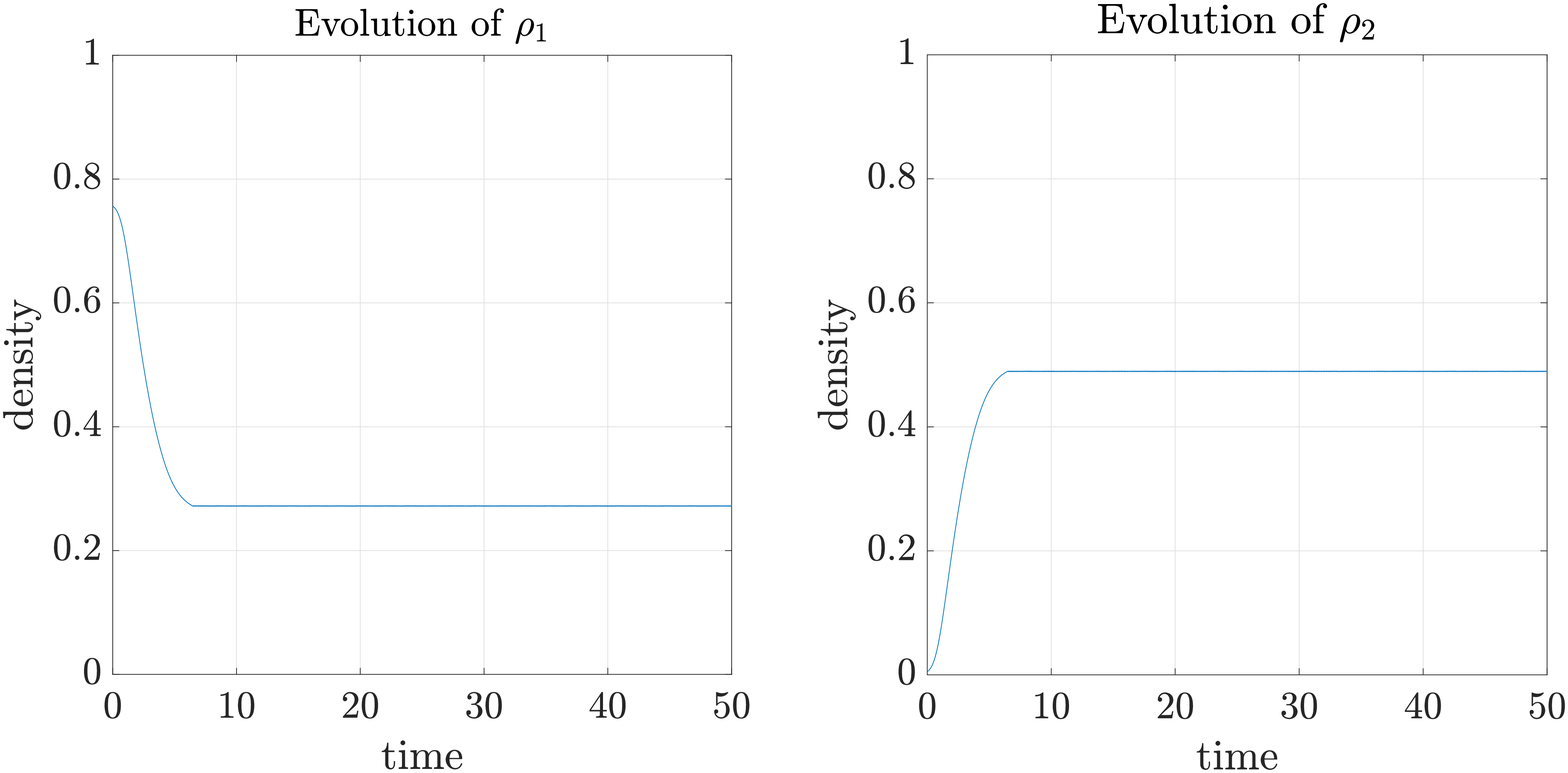} \includegraphics[width=0.49\textwidth]{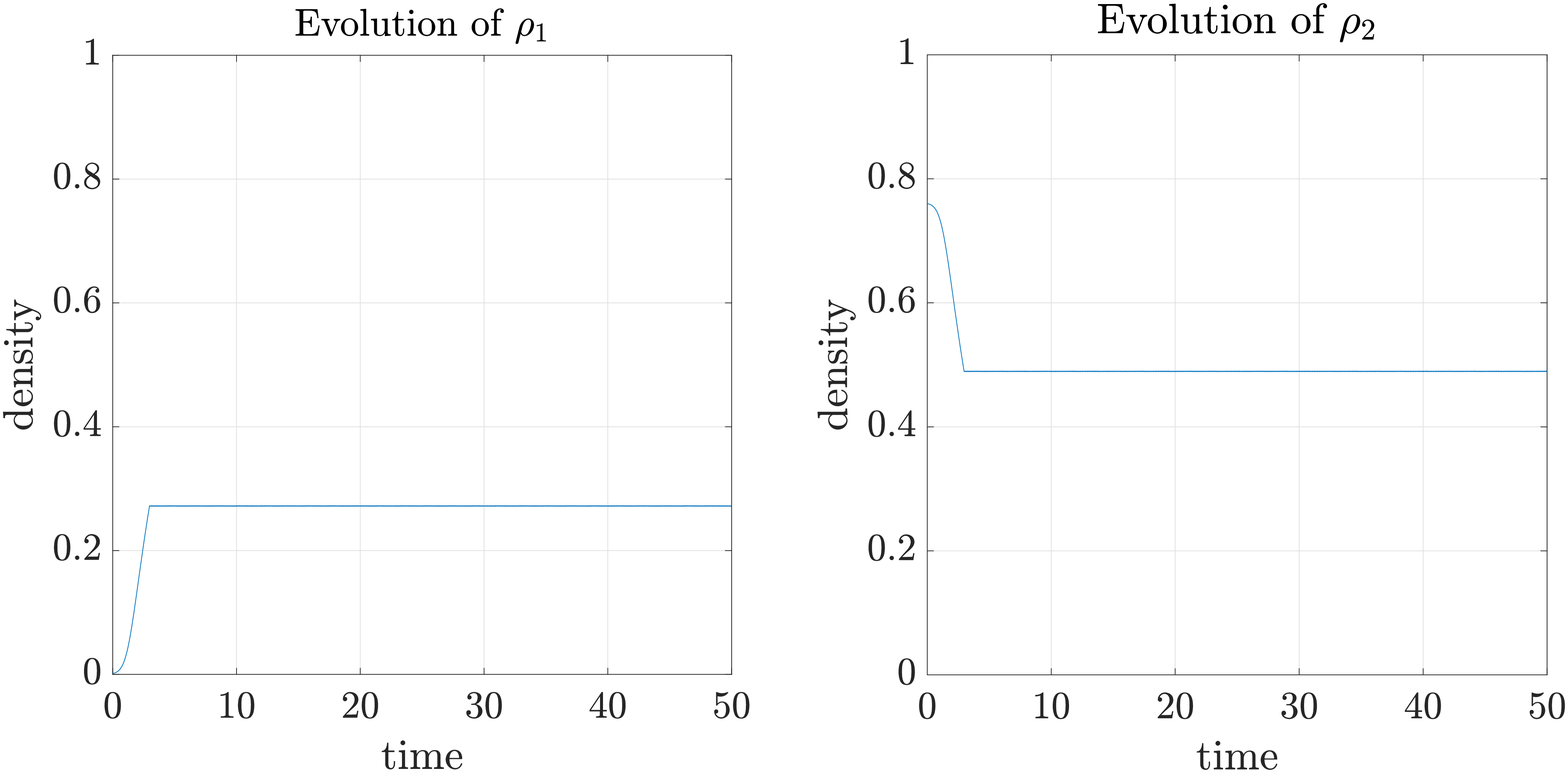}
\caption[Global in space perturbation of the equilibrium solution.]{Global in space perturbation of the equilibrium solution $(\overline{\rho}_1,\overline{\rho}_2)=(0.27,0.49)$ of type (A). Left: the perturbation is $(\rho_1,\rho_2)=(\overline{\rho}_1+\epi\rho_0,\overline{\rho}_2-\epi\rho_0)$, with $\epi\rho_0=0.485$. Right: the perturbation is $(\rho_1,\rho_2)=(\overline{\rho}_1+\epi\rho_0,\overline{\rho}_2-\epi\rho_0)$, with $\epi\rho_0=-0.27$.}\label{evoA}
\end{figure} 

Here, the space domain is $[-0.5,0.5]$ with periodic boundary conditions, the final time is $T=50$, and $\Delta x =0.003$. Clearly the density remains uniform in space. In Fig.~\ref{evoA} we show the evolution in time of  the densities $\rho_1$ and $\rho_2$ on the road length. We observe that $ \rho_j(t) \to \overline{\rho}_j$, $j=1,2$, as $t\to+\infty$. Therefore, $\rho_j(t) \to \overline{\rho}_j$, $j=1,2$, as $t\to\infty$ shows the convergence of the system to the original steady state. This confirms the theoretical finding on the global asymptotic stability of equilibria of type (A).

\subsection{Local  perturbation in space of equilibrium solutions}

In the following test, we study the effects of a local  perturbation in space. In particular, the perturbation of the equilibrium state $(\overline{\rho}_1,\overline{\rho}_2)$ is a bump in the density located around $x=0$ and it is modeled with a Gaussian profile. More precisely, we consider the following initial datum:

\begin{equation} \label{eq:local:pert}
\begin{split}
&(\rho_1(x,0),\rho_2(x,0))=\left(\overline{\rho}_1+G(x),\overline{\rho}_2-G(x) \right) \\ &\text{with } (\overline{\rho}_1,\overline{\rho}_2)=(0.142,0.400),
\end{split}
\end{equation}
where
\begin{equation} \label{eq:gaussian} 
	G(x)=\exp(-100x^2)\cdot R \quad \text{with } R=\begin{cases} \min(\overline{\rho}_1, 1-\overline{\rho}_2) & \mbox{if } \overline{\rho}_1\geqslant \overline{\rho}_2, \\ \min(\overline{\rho}_2, 1-\overline{\rho}_1)  & \mbox{otherwise.}  \end{cases}
\end{equation}

The equilibrium $(\overline{\rho}_1,\overline{\rho}_2)$ is of type (A), and precisely it is  obtained with the equilibrium speed $v^{\text{eq}}=0.6$. The space domain is again $[-0.5,0.5]$ with periodic boundary conditions, the final time is $T=5$, $\Delta x = 0.01$, and $CFL=0.9$.

\begin{figure}[ht]
\centering
\subfigure[Density profiles at different snapshots: $t=0$, $t=0.2024$, $t=0.5088$, $t=1.0035$, $t=2.5110$, $t=T=5$. Red: $\rho_1$, blue: $\rho_2$.]{
\includegraphics[width=0.9\textwidth]{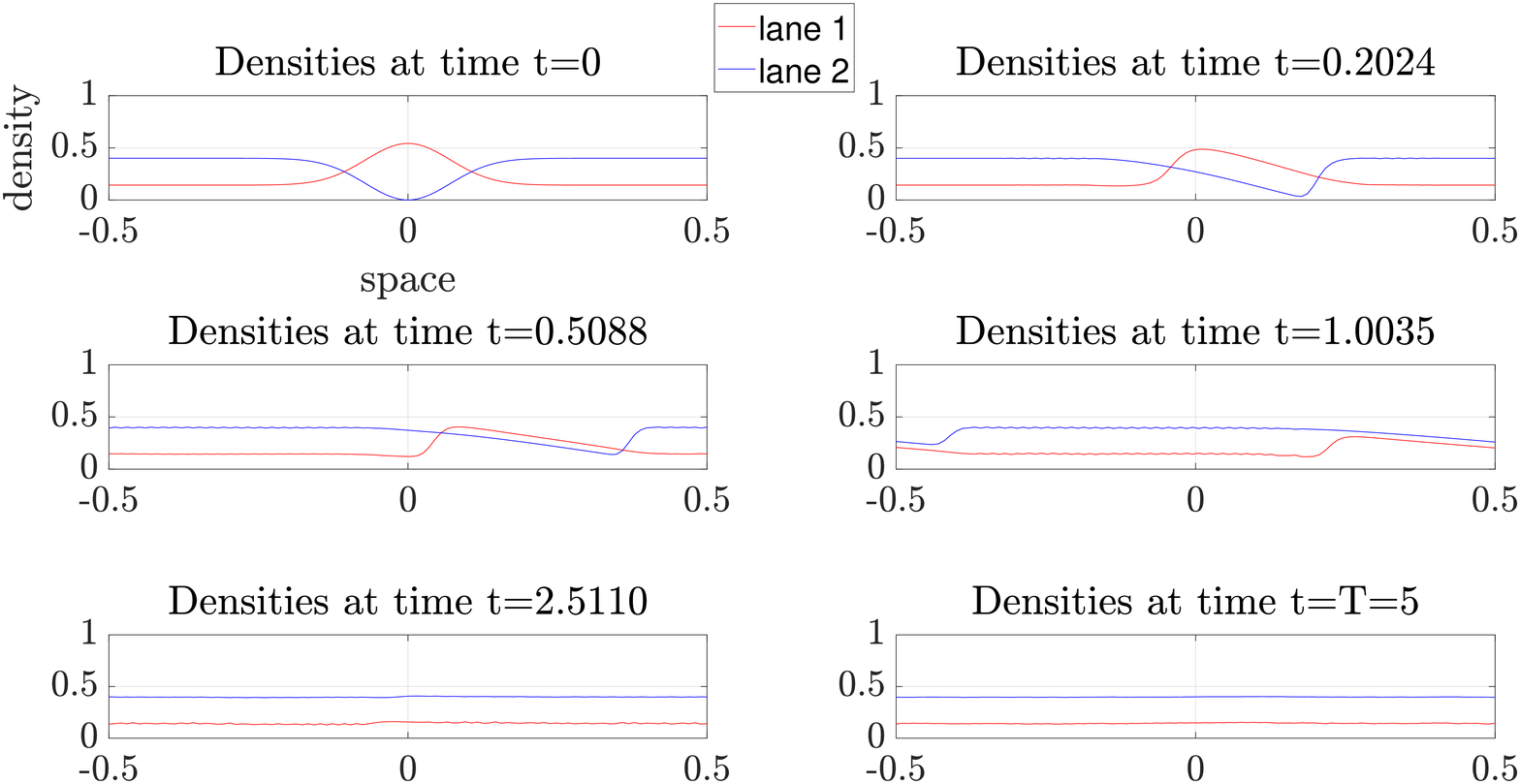}}
\subfigure[Evolution in time of the density profiles in lane 2 (left panel) and lane 1 (right panel).]{
\includegraphics[width=0.9\textwidth]{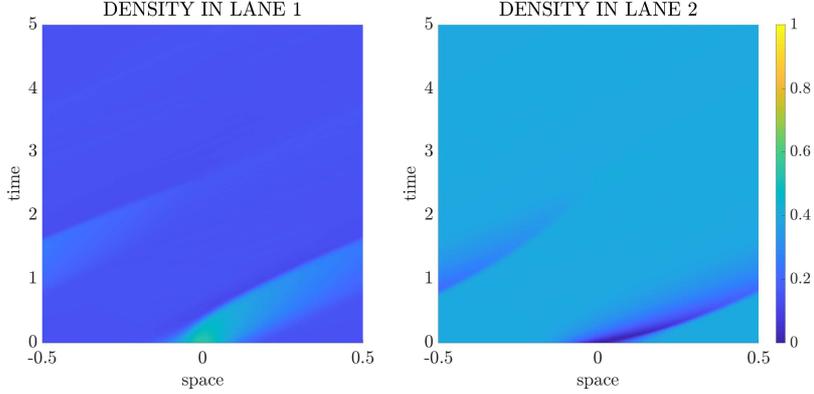}}
\caption[Local in space perturbation of the equilibrium solution.]{Local in space perturbation of the equilibrium solution with model \eqref{modello}. }\label{gauss}
\end{figure} 

The density profiles at different snapshots on both lanes are shown in Fig.~\ref{gauss}(a), whereas the evolution in time can be observed in Fig.~\ref{gauss}(b). Computing the mean value of the densities $\rho_1$ and $\rho_2$,  as
\begin{equation} \label{eq:mean:density}
	\langle \rho_j \rangle(t) \approx \int_{-\frac12}^{\frac12} \rho_j(x,t)\, \mathrm{d}x, \quad j=1,2,
\end{equation} we obtain that $\langle \rho_1 \rangle(T)=0.144$, with standard deviation
$4\cdot10^{-3}$, and $\langle \rho_2 \rangle(T)=0.399$ with standard deviation $2\cdot10^{-3}$. We conclude that the numerical results show the convergence of the system to the original equilibrium $(\overline{\rho}_1,\overline{\rho}_2)$. Therefore, also in the case of a local perturbation, the system behaves as studied in the case of a homogeneous perturbation around the steady states.

\subsection{Lane closure}
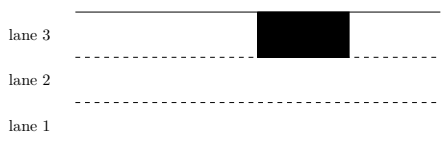
\begin{figure}[H]
\centering
\scalebox{0.6}{
\begin{tikzpicture}
\draw (0,3)--(8,3);
\draw [dashed](0,2)--(8,2);
\draw [dashed](0,1)--(8,1);
\draw (0,0)--(8,0);
\draw[fill=black] (4,2) rectangle ++(2,1);
\path(-1,0.5)node{lane $1$};
\path (-1,1.5)node{lane $2$};
\path (-1,2.5)node{lane $3$};
\end{tikzpicture}}
\caption{Schematic for traffic flow through local closure of lane 3.}\label{chiusura}
\end{figure}

Here we present an example of a lane closure on a three-lane road. In particular a stretch  of the fastest lane, i.e.~lane 3, is closed due to, e.g., an accident or the presence of a work-in-progress area. The setup for this simulation is illustrated in Fig.~\ref{chiusura}. On the three-lane road we prescribe different desired velocity profiles such that $v_3(\rho)>v_2(\rho)>v_1(\rho)$, for $\rho\in[0,\rho^{\text{max}})$ and $v_3(\rho^{\text{max}})=v_2(\rho^{\text{max}})=v_1(\rho^{\text{max}})=0$.

We consider a road of unit length modeled by the space domain $I=[-0.5,0.5]$ with closure in lane 3 located in $C=[0,0.25]$. Linear velocities as in the previous tests are prescribed with $v_3^{\max}=1$, $v_2^{\max}=0.7$, and $v_1^{\max}=0.6$. We present two tests with uniform initial densities given by
\begin{equation}
\begin{array}{l}
\rho_1(x,0)=0.4,\,\, x \in I\,\, \mbox{(test 1 \& 2),} \\ 
\rho_2(x,0)=0.6, \,\,x\in I\,\, \mbox{(test 1);}\,\,\,\, \rho_2(x,0)=0.1, \,\,x\in I\,\, \mbox{(test 2),} \\
\rho_3(x,0)=\begin{cases} 0.2 & x<0 \\ 0 & x>0.25. \end{cases} \,\, \mbox{(test 1 \& 2).}
\end{array}
\end{equation}
Moreover, we assume free flow boundary conditions on $x=0.5$ and Dirichlet boundary conditions on $x=-0.5$ given, for $j=1,2,3$, by
\begin{equation*}
    \rho_j(-0.5,t)=\rho_j(-0.5,0) \quad\forall t>0.
\end{equation*}
We evolve the model~\eqref{modello} with $J=3$, up to final time $T=1.2$. We consider $\Delta x=0.001$ and $CFL=0.99$.

\begin{figure}[ht]
\centering
\subfigure[Test 1.]{
\includegraphics[width=0.9\textwidth]{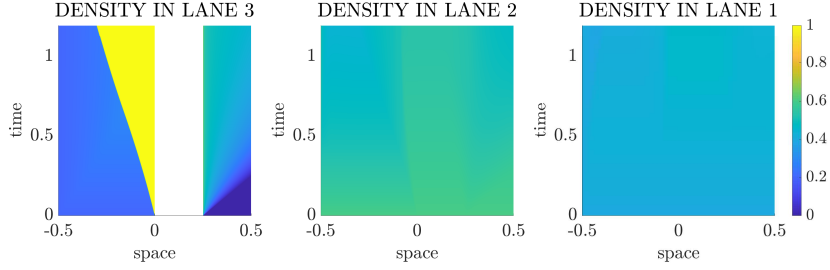}}
\subfigure[Test 2.]{
\includegraphics[width=0.9\textwidth]{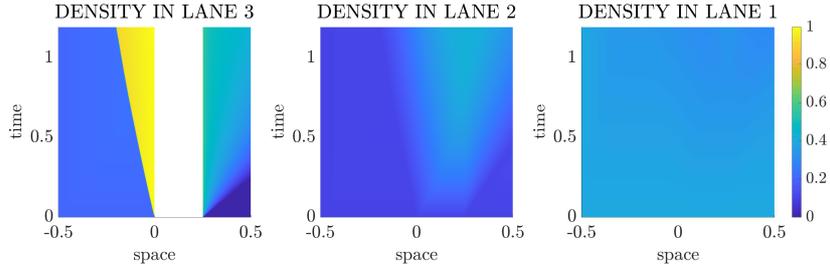}}
\caption{Lane closure numerical test.}\label{2d}
\end{figure}

Fig.~\ref{2d} shows the evolution of the densities in the three lanes. We observe that the local closure of lane 3 leads to an increase of the densities in the two rightmost lanes that can also create bottlenecks depending on the initial congestion of the road. Both tests show the presence of a queue in lane 3 that evolves backward in space. In test 1 (the congested case in lane 2) vehicles migrate to lane 2 where the density tends to the critical value $\mu=0.5$,  while in test 2 (the uncongested case in lane 2) the density in lane 2 tends to $\mu$ only  locally near the beginning of the closure due to the different inflow condition. Furthermore, we note that the  density in lane 1, the slowest one, remains close to the value $\mu$ because vehicles  in lane 2 cannot jump in lane 1.

\section{Conclusions} \label{sec:conclusion}
In this paper, we have presented a new first-order multi-lane macroscopic model for traffic flow, which has been derived as continuum limit of a microscopic follow-the-leader framework where we have proposed specific lane changing conditions. Exploiting the evolution of the discrete densities, we have described the effects of lane changes in terms of macroscopic quantities and without postulating specific scalings in space and time for the derivation of the final system of partial differential equations. The system thus obtained is a system of balance laws with source terms that refer to the lane changing and that model the increase or decrease of the density in a given lane, in terms of itself. After this derivation a detailed stability analysis of the system equilibria was performed. The considered equilibria refer to a situation where traffic flow is constant in all lanes and there are no lane changes. We have studied the stability of these equilibria under global, i.e.~space homogeneous, perturbations proving the existence of asymptotically stable, globally asymptotically stable and marginally stable equilibria. We have performed some numerical tests devoted to the comparison between the microscopic model and the corresponding multi-lane macroscopic limit. Furthermore, we have studied numerically the evolution of perturbed equilibria using global perturbation, confirming the results of the theoretical stability analysis. We have also shown that the steady states behave similarly under local perturbations. Finally, we have proposed a numerical test aimed at studying the behavior of the model in a real situation such as the closure of part of a lane on a three-lane road. The evolution of the densities show realistic flow of vehicles with the lane changes dynamics introduced here.

\section{Appendix: proof of Theorem \ref{teorema}}
Let us illustrate the analysis of each equilibrium class.

\textbf{Steady state of type (A).} In this case, the monotonicity of the velocity functions implies that
\begin{equation*}\label{condA}
v_1^\mu < v^{\text{eq}} \leqslant v_1^{\max}. 
\end{equation*}
We consider two subcases depending on the sign of $\epi$:
\begin{itemize}
\item[(A1)] If $\epi>0$ then:
\begin{equation*}
\begin{aligned}
\rho_1>\overline{\rho}_1 \quad &\text{and} \quad v_1(\rho_1)<v^{\text{eq}}, \quad 
\rho_2<\mu \quad &\text{and} \quad v_2(\rho_2)>v^{\text{eq}}.
\end{aligned}
\end{equation*}
Thus, lane changes from lane 2 to lane 1 do not occur, i.e.~$\pi^{2 \to 1}(\rho_2,\rho_1)=0$, because the incentive condition is not satisfied.
\item[(A2)] Similarly in the case $\epi<0$, in which lane changes from lane 1 to lane 2 do not occur.
\end{itemize}

Let us study in detail  case (A1). 
We have that the evolution of the perturbed states $(\rho_1,\rho_2)$ follows the following system of balance laws:
\begin{equation*}
\begin{cases}
\de_t \rho_1 + \de_x f_1(\rho_1) =  - \pi^{1 \to 2}(\rho_1,\rho_2) A(\rho_1,\rho_2)\rho_2, \\
\de_t \rho_2 + \de_x f_2(\rho_2) = \pi^{1 \to 2}(\rho_1,\rho_2) A(\rho_1,\rho_2)\rho_2,
\end{cases}
\end{equation*}
and substituting~\eqref{perturbazione_modello} we obtain
\begin{equation*}
    \resizebox{1\textwidth}{!}{%
        $\begin{cases} \de_t (\overline{\rho}_1+\epi \rho) + \de_x f_1(\overline{\rho}_1+\epi \rho) =  - \pi^{1 \to 2}(\overline{\rho}_1+\epi \rho,\overline{\rho}_2-\epi \rho) A(\overline{\rho}_1+\epi \rho,\overline{\rho}_2-\epi \rho)(\overline{\rho}_2-\epi \rho),\\ \de_t (\overline{\rho}_2-\epi \rho) + \de_x f_2(\overline{\rho}_2-\epi \rho) =+\pi^{1 \to 2}(\overline{\rho}_1+\epi \rho,\overline{\rho}_2-\epi \rho) A(\overline{\rho}_1+\epi \rho,\overline{\rho}_2-\epi \rho)(\overline{\rho}_2-\epi \rho). \end{cases}$      
        }
     \end{equation*} 
If we linearize the previous equations around the steady state $(\overline{\rho}_1,\overline{\rho}_2)$, namely
\begin{align*}
\pi^{1 \to 2}(\overline{\rho}_1+\epi \rho,\overline{\rho}_2-\epi \rho) &= \pi^{1 \to 2}(\overline{\rho}_1,\overline{\rho}_2)+\epi \nabla  \pi^{1 \to 2}(\overline{\rho}_1,\overline{\rho}_2) \cdot (\rho,-\rho)^T +O(\epi^2) \\ &= -\epi g'(\overline{\rho}_2)\rho + O(\epi^2), \\
A(\overline{\rho}_1+\epi \rho,\overline{\rho}_2-\epi \rho) &= A(\overline{\rho}_1,\overline{\rho}_2)+\epi \nabla A(\overline{\rho}_1,\overline{\rho}_2)\cdot (\rho,-\rho)^T +O(\epi^2) \\ &=A(\overline{\rho}_1,\overline{\rho}_2)+\epi(\de_{\rho_1}A(\overline{\rho}_1,\overline{\rho}_2)-\de_{\rho_2}A(\overline{\rho}_1,\overline{\rho}_2))\rho+O(\epi^2),
\end{align*}
we get the evolution of the perturbation $\rho=\rho(t)$ as
\begin{equation*}
\de_t \rho =  \overline{\rho}_2 A(\overline{\rho}_1, \overline{\rho}_2)g'( \overline{\rho}_2)\rho + O(\epi)
\end{equation*}
with analytical solution
\begin{equation} \label{eq:exp:decay:rho}
\rho(t)=\rho_0 \exp\left(\overline{\rho}_2 A(\overline{\rho}_1, \overline{\rho}_2)g'( \overline{\rho}_2) t\right).
\end{equation}
Recalling that $g'(r)<0 \ \forall r$, we find that $\rho(t)\to 0$ as $t\to +\infty$, and thus $(\rho_1,\rho_2)\to(\overline{\rho}_1,\overline{\rho}_2)$ as $t\to +\infty$. In other words, the perturbed state converges to the steady state of type (A1) under the perturbation~\eqref{perturbazione_modello}.

With similar computations we can obtain the same result for  subcase (A2), therefore we omit explicit presentation.

\textbf{Steady state of type (B).} In this case it is possible to distinguish between various situations due to the monotonicity assumption of the velocity functions.
\begin{itemize}
\item[(B1)] If $\overline{\rho}_1\geqslant \mu$, \eqref{eq:speed:functions} implies that also $\overline{\rho}_2>\mu$ and $\overline{\rho}_2>\overline{\rho}_1$; then, in this case
\begin{equation*}
0 \leqslant v^{\text{eq}} \leqslant v_1^{\mu}.
\end{equation*}
\item[(B2)] If $\overline{\rho}_2\geqslant \mu$, instead, the only interesting case is provided when  $\mu \leqslant \overline{\rho}_2<\rho_2^\mu$. In fact, in this case,  \eqref{eq:speed:functions} implies that $\overline{\rho}_1<\mu$ and thus
\begin{equation}\label{defb2}
v_1^\mu < v^{\text{eq}} \leqslant v_2^{\mu}.
\end{equation}
On the other hand if $\overline{\rho}_2\geqslant \rho_2^\mu$ then $\overline{\rho}_1\geqslant \mu$ and this case reduces to (B1).
\end{itemize}

As in the steady state of type (A), we study the stability under a perturbation of the form~\eqref{perturbazione_modello}. 

Let us analyze first  case (B1). Given an equilibrium $(\overline{\rho}_1,\overline{\rho}_2)$ of type (B1), we can distinguish the following sub-cases depending on the sign of its perturbation:
\begin{itemize}
\item In the case $\epi>0$:
\begin{itemize}
\item[(I)] if the perturbed densities $(\rho_1,\rho_2)$ are both above the critical value, i.e.~$\mu < \rho_2 < \rho_1$, then $(\rho_1,\rho_2)$ is an equilibrium of type (C) for the system because lane changing is not possible due to the safety criterion;
\item[(II)] if, instead, the perturbed density in the second lane is less than the critical value, i.e.~$\rho_2<\mu$, then the system converges to a different equilibrium solution $(\tilde{\rho}_1,\tilde{\rho}_2)$, which is of type (C) with $\tilde{\rho}_2=\mu$;
\end{itemize}
\item Similar results hold for the case $\epi<0$.
\end{itemize}
We prove the previous statements in the following.

Case (I) implies that lane changes are not possible either from lane 1 to lane 2 or from lane 2 to lane 1, because the safety criterion is violated. Thus, $\pi^{1 \to 2}(\rho_1,\rho_2)=\pi^{2\to 1}(\rho_2,\rho_1)=0$  $ \forall(x,t)$, and it is easy to check that the evolution of the perturbation $\rho$ satisfies
\begin{equation*}
\epi \de_t \rho =  0  \,\,\mbox{   which implies   }\,\, \rho(t) = \rho_0.
\end{equation*}
Then, $(\rho_1,\rho_2)=(\overline{\rho}_1+\epi\rho_0,\overline{\rho}_2-\epi\rho_0)$ is an equilibrium of type (C) because the densities are both above the critical value $\mu$ and the speeds corresponding to the new states are different. 

Case (II) implies that lane changes are possible only from lane 1 to lane 2, because of the safety condition. Thus, $\pi^{1 \to 2}(\rho_1,\rho_2)>0$ and $\pi^{2\to 1}(\rho_2,\rho_1)=0$. Repeating the same computations as for the steady states of type (A), it is possible to show that we obtain again the exponential decay~\eqref{eq:exp:decay:rho} of the perturbation $\rho$. However, in this case the density $\rho_1$ decreases as long as $\pi^{1 \to 2}(\rho_1,\rho_2)>0$. 
Therefore, there exists a time $\overline{t}$ such that $\pi^{1\to 2}\to0$ as $t\to\overline{t}$, and
\begin{equation*}
\rho(t) =
\begin{cases}
\rho_0 \exp\left(\overline{\rho}_2 A(\overline{\rho}_1, \overline{\rho}_2)g'( \overline{\rho}_2) t\right), & t<\overline{t} \\[1ex]
\frac{1}{\epi}(\overline{\rho}_2-\mu), & t\geqslant \overline{t}
\end{cases}
\end{equation*}
where the behavior for $t\geqslant \overline{t}$ is obtained by imposing the mass conservation principle.

Now, we focus on the case (B2). Given an equilibrium $(\overline{\rho}_1,\overline{\rho}_2)$ of type (B2), we can distinguish the following subcases depending on the sign of its perturbation:
\begin{itemize}
\item In the case $\epi>0$ we have: 
\begin{itemize}
\item[(I)] if the perturbed state $(\rho_1,\rho_2)$ is such that $\rho_1<\mu$ and $\rho_2>\mu$, then the system converges to a different steady state $(\tilde{\rho}_1,\tilde{\rho}_2)$, which is of type (D);
\item[(II)] if the perturbed state $(\rho_1,\rho_2)$ is such that $\rho_1>\mu$ and $\rho_2>\mu$, then the system converges to a different steady state $(\tilde{\rho}_1,\tilde{\rho}_2)$, which is of type (C);
\item[(III)] if the perturbed state $(\rho_1,\rho_2)$ is such that $\rho_2 \leqslant \mu$, then the system converges to a different steady state $(\tilde{\rho}_1,\tilde{\rho}_2)$, which is of type (C) with $\tilde{\rho}_2=\mu$;
\end{itemize}
\item in the case $\epi<0$: 
\begin{itemize}
\item[(IV)] independently of the initial perturbed state $(\rho_1,\rho_2)$, the system returns to the equilibrium $(\overline{\rho}_1,\overline{\rho}_2)$ of type (B2). 
\end{itemize}
\end{itemize}

We observe that in case (I) lane changes are not possible, either from lane 1 to lane 2 because of the safety criterion, or from lane 2 to lane 1 because of the incentive criterion. Therefore, the perturbation $\rho$ remains constant in time.

Case (II) is similar, since there are no lane changes because of the safety criterion. 

In case (III), instead, only lane changes from lane 1 to lane 2 are possible and until the safety criterion is verified, i.e.~until $\rho_2$ remains smaller than the critical value $\mu$. Therefore, the analysis of this case is similar to the analysis of case (II) for a steady state of type (B1).

Finally, for case (IV) we note that, since $\rho_1<\mu$ and $\rho_2>\mu$, only lane changing from lane 2 to lane 1 can occur. The velocity in lane 2 increases until its value becomes equal to the velocity in lane 1 which is decreasing. In this case, the system converges to the original equilibrium $(\overline{\rho}_1, \overline{\rho}_2)$.

\textbf{Steady state of type (C).} If the perturbation~\eqref{perturbazione_modello} of an equilibrium state $(\overline{\rho}_1, \overline{\rho}_2)$ of type (C) is characterized by a small $\epi \rho$ such that $\rho_1>\mu$ and $\rho_2>\mu$, then $\partial_t \rho = 0$ and $(\rho_1,\rho_2)=(\overline{\rho}_1+\epi\rho,\overline{\rho}_2-\epi\rho)$ establishes an equilibrium of type (C) for the system. Otherwise depending on the size of the initial perturbation $\epi \rho_0$ the system converges  to a different equilibrium either of type (C) characterized by the density on one lane equal to $\mu$, or of type (D), or of type (B1), or of type (B2).   

\textbf{Steady state of type (D).} We observe that an equilibrium state $(\overline{\rho}_1, \overline{\rho}_2)$ of type (D) is characterized by $\overline{\rho}_1<\mu$, so that the safety condition for a lane change from lane $2$ to lane $1$ is satisfied. However, the lane change cannot occur because at the same time  $v_1(\overline{\rho}_1)<v_2(\overline{\rho}_2)$ holds. 
This situation is  possible only if
\begin{equation}\label{caseD}
\mu \leqslant \overline{\rho}_2 \leqslant \rho_2^\mu \qquad \text{and} \qquad (v_1)^{-1}(v_{2}(\overline{\rho}_2))\leqslant\overline{\rho}_1 \leqslant \mu.
\end{equation}
If one has that $\overline{\rho}_1=(v_1)^{-1}(v_{2}(\overline{\rho}_2))$, case (D) reduces to case (B2).
On the contrary, if $\overline{\rho}_1\neq (v_1)^{-1}(v_{2}(\overline{\rho}_2))$ and the perturbed state $(\rho_1,\rho_2)$ (cf.~\eqref{perturbazione_modello}) satisfies~\eqref{caseD}, then $\partial_t \rho = 0$ because the lane changing conditions are not satisfied, and $(\rho_1,\rho_2)$ establishes a new equilibrium of type (D) for the system. Otherwise, if $(\rho_1,\rho_2)$ does not satisfy~\eqref{caseD}, depending on the size of the perturbation the system converges to an equilibrium state either of type (C), or of type (B2).

\textbf{Steady state of type (E).} In this case the equilibrium state is characterized by a zero density in the slowest lane, whereas the density in the fastest lane takes values in $[0,v_2^{-1}(v_1^{\text{max}})]$. Perturbation~\eqref{perturbazione_modello} is admissible only for $\epi>0$, thus defining a perturbed state $(\rho_1,\rho_2)$ of the system which makes possible only lane changes from lane 1 to lane 2. With the same arguments of case (A1), we can prove that the perturbation $\rho(t)$ decays exponentially and $(\rho_1,\rho_2)\to(\overline{\rho}_1,\overline{\rho}_2)$ as $t\to\infty$. This result explains mathematically the situation in which all vehicles migrate from the slowest lane to the fastest one. \hfill $\square$

\section*{Funding and acknowledgments}
{\footnotesize M.H. thanks the Deutsche Forschungsgemeinschaft (DFG, German Research Foundation) for the financial support through 320021702/GRK2326, 333849990/IRTG-2379, B04, B05 and B06 of CRC1481, HE5386/18-1,19-2,22-1,23-1, ERS SFDdM035 and under Germany’s Excellence Strategy EXC-2023 Internet of Production 390621612 and under the Excellence Strategy of the Federal Government and the Länder. M.H. and G.P. acknowledge support through the EU DATAHYKING project. M.P. and G.P. acknowledge support through Ateneo Sapienza project 2019 ``Metodi numerici per problemi evolutivi, networks ed applicazioni'', 2020 ``Algoritmi e modelli per sistemi di natura iperbolica, networks e applicazioni'', and 2021 ``Evolutionary problems: analysis techniques and construction of numerical solutions''. This work was also carried out within the MUR (Ministry of University and Research) PRIN-2017 project ``Innovative Numerical Methods for Evolutionary Partial Differential Equations and Applications'' (number 2017KKJP4X). M.P. and G.P. are members of the INdAM Research Group GNCS. M. P. wishes to thank: Giuseppe Visconti for the valuable help and fruitful discussions; and Elisa Iacomini for the  useful suggestions during the stay at the RWTH Aachen University.}

\printbibliography

@ARTICLE{4806347,
  author={{Goebel}, R. and {Sanfelice}, R. G. and {Teel}, A. R.},
  journal={IEEE Control Systems Magazine}, 
  title={Hybrid dynamical systems}, 
  year={2009},
  volume={29},
  number={2},
  pages={28-93},}

@inproceedings{piccoli1998hybrid,
  title={Hybrid systems and optimal control},
  author={Piccoli, B.},
  booktitle={Proceedings of the 37th IEEE Conference on Decision and Control (Cat. No. 98CH36171)},
  volume={1},
  pages={13-18},
  year={1998},
  organization={IEEE}
}

@article{garavello2005hybrid,
  title={Hybrid necessary principle},
  author={Garavello, M. and Piccoli, B.},
  journal={SIAM J. Control Optim.},
  volume={43},
  number={5},
  pages={1867-1887},
  year={2005},
  publisher={SIAM}
}

@article{PiuPuppo2022,
title = {Stability analysis of microscopic models for traffic flow with lane changing},
journal = {Netw. Heterog. Media},
volume = {17},
number = {4},
pages = {495-518},
year = {2022},
issn = {1556-1801},
author = {Piu, M. and Puppo, G.}
}

@Article{SukhinovaTrapeznikovaChetverushkinChurbanova2009,
	author = {Sukhinova, A.B. and Trapeznikova, M.A. and Chetverushkin, B. N. and Churbanova, N. G.},
	title = {Two-{D}imensional {M}acroscopic {M}odel of {T}raffic {F}lows},
	Journal = {Mathematical Models and Computer Simulations},
	Volume = {1},
	Number = {6},
	Pages = {669--676},
	Year = {2009}
}

@article{ChiriGongPiccoli,
    author = {Chiri, M. T. and Gong, X. and Piccoli, B.},
    title = {Mean-field limit of a hybrid system for multi-lane car-truck traffic},
    year = {2022},
    journal = {Netw. Heterog. Media},
    %note = {To appear}
}

@Unpublished{GongPiccoliVisconti,
	author = {Gong, X. and Piccoli, B. and Visconti, G.},
	title = {Mean-field limit of a hybrid system for multi-lane multi-class traffic},
	year = {2022},
	%note = {Preprint arXiv:2007.14655}
}

@article{SongKarni2019,
author = {Song, J. and Karni, S.},
year = {2019},
title = {A Second Order Traffic Flow Model with Lane Changing},
journal = {J. Sci. Comput.},
pages = {1429-1445},
volume = {81},
issue = {3}
}

@article{HoldenRisebro2019,
author = {Holden, H. and Risebro, N. H.},
title = {Models for Dense Multilane Vehicular Traffic},
journal = {SIAM J. Math. Anal.},
volume = {51},
number = {5},
pages = {3694-3713},
year = {2019}
}

@article{GoatinRossi2019,
author = {Goatin, P. and Rossi, E.},
title = {A MultiLane Macroscopic Traffic Flow Model for Simple Networks},
journal = {SIAM J. Appl. Math.},
volume = {79},
number = {5},
pages = {1967-1989},
year = {2019}
}

@article{albi2019vehicular,
	title={Vehicular traffic, crowds, and swarms: From kinetic theory and multiscale methods to applications and research perspectives},
	author={Albi, G and Bellomo, N and Fermo, L and Ha, S-Y and Kim, J and Pareschi, L and Poyato, D and Soler, J},
	journal={Mathematical Models and Methods in Applied Sciences},
	volume={29},
	number={10},
	pages={1901--2005},
	year={2019},
	publisher={World Scientific}
}

@Article{aw2000SIAP,
	Title                    = {Resurrection of ``second order'' models of traffic flow},
	Author                   = {Aw, A. and Rascle, M.},
	Journal                  = {SIAM J. Appl. Math.},
	Year                     = {2000},
	Number                   = {3},
	Pages                    = {916-938 (electronic)},
	Volume                   = {60},
	
	Fjournal                 = {SIAM Journal on Applied Mathematics},
	ISSN                     = {0036-1399},
	Mrclass                  = {35L65 (90B06)},
	Mrnumber                 = {1750085 (2001a:35111)},
	Mrreviewer               = {Yu. G. Rykov}
}

@Article{aw2002SIAP,
	Title                    = {Derivation of continuum traffic flow models from microscopic follow-the-leader models},
	Author                   = {Aw, A. and Klar, A. and Materne, T. and Rascle, M.},
	Journal                  = {SIAM J. Appl. Math.},
	Year                     = {2002},
	Number                   = {1},
	Pages                    = {259-278},
	Volume                   = {63},
	
	Coden                    = {SMJMAP},
	Fjournal                 = {SIAM Journal on Applied Mathematics},
	ISSN                     = {0036-1399},
	Mrclass                  = {35L45 (35L65 90B20)},
	Mrnumber                 = {1952895 (2003m:35148)},
	Mrreviewer               = {Alan Jeffrey}
}

@book {MR1600250,
    AUTHOR = {Haberman, R.},
     TITLE = {Mathematical models},
    SERIES = {Classics in Applied Mathematics},
    VOLUME = {21},
      %note = {Mechanical vibrations, population dynamics, and traffic flow, An introduction to applied mathematics, Reprint of the 1977 original},
 PUBLISHER = {Society for Industrial and Applied Mathematics (SIAM), Philadelphia, PA},
      YEAR = {1998},
     PAGES = {xviii+402},
   MRCLASS = {00A69 (00A71 34-01 70-01 90-01)},
  MRNUMBER = {1600250}
}

@ARTICLE{Helbing20011067,
author={Helbing, D.},
title={Traffic and related self-driven many-particle systems},
journal={Reviews of Modern Physics},
year={2001},
volume={73},
number={4},
pages={1067-1141}
}

@Article{HertyMoutariVisconti2018,
	Title                    = {Macroscopic modeling of multilane motorways using a two-dimensional second-order model of traffic flow},
	Author                   = {Herty, M. and Moutari, S. and Visconti, G.},
	Journal                  = {SIAM J. Appl. Math.},
	Year                     = {2018},
	Number                   = {4},
	Pages                    = {2252-2278},
	Volume                   = {78}
}

@Article{lighthill1955PRSL,
	Title                    = {On kinematic waves. {II}. {A} theory of traffic flow on long crowded roads},
	Author                   = {Lighthill, M. J. and Whitham, G. B.},
	Journal                  = {Proc. Roy. Soc. London. Ser. A.},
	Year                     = {1955},
	Pages                    = {317-345},
	Volume                   = {229},
	
	Fjournal                 = {Proceedings of the Royal Society. London. Series A. Mathematical, Physical and Engineering Sciences}
}

@Article{richards1956OR,
	Title                    = {Shock waves on the highway},
	Author                   = {Richards, P. I.},
	Journal                  = {Operations Res.},
	Year                     = {1956},
	Pages                    = {42-51},
	Volume                   = {4}
}

@InCollection{piccoli2009ENCYCLOPEDIA,
	Title                    = {Vehicular traffic: {A} review of continuum mathematical models},
	Author                   = {Piccoli, B. and Tosin, A.},
	Booktitle                = {Encyclopedia of Complexity and Systems Science},
	Publisher                = {Springer},
	Year                     = {2009},
	
	Address                  = {New York},
	Editor                   = {Meyers, R. A.},
	Pages                    = {9727-9749},
	Volume                   = {22}
}

@Article{Zhang2002,
	Title                    = {A non-equilibrium traffic model devoid of gas-like behavior},
	Author                   = {Zhang, H. M.},
	Journal                  = {Transport. Res. B-Meth.},
	Year                     = {2002},
	Number                   = {3},
	Pages                    = {275-290},
	Volume                   = {36}
}

\end{document}